\documentclass[3p, sort&compress, times]{elsarticle}
\nopreprintlinetrue

\usepackage[utf8]{inputenc}
\usepackage{amssymb, amsmath, amsthm, mathtools, upgreek}
\usepackage{array,booktabs}
\newcolumntype{C}{>{$}c<{$}}
\newcolumntype{R}{>{$}r<{$}}
\usepackage{algpseudocodex,algorithm}
\usepackage{xcolor}
\definecolor{tumblue}{rgb}{0,0.3961,0.7412}
\colorlet{blue}{tumblue}
\usepackage{hyperref}
\usepackage[capitalise]{cleveref}
\hypersetup{
    pdftitle={A lattice Boltzmann method for Biot's consolidation model of linear poroelasticity},
    pdfauthor={S. B. Lunowa, B. Wohlmuth},
    colorlinks=true
}

\newcommand*{\N}{\mathbb{N}}
\renewcommand{\vec}[1]{\boldsymbol{#1}}
\newcommand{\tensor}[1]{\boldsymbol{\mathbf{#1}}}
\providecommand\given{}
\newcommand\SetSymbol[1][]{\nonscript\:#1\vert\allowbreak\nonscript\:\mathopen{}}
\DeclarePairedDelimiterX{\Set}[1]{\{}{\}}{\renewcommand{\given}{\SetSymbol[\delimsize]}#1}
\DeclarePairedDelimiterX{\Norm}[1]{\lVert}{\rVert}{#1}

\begin{document}

\begin{frontmatter}

\title{A lattice Boltzmann method for Biot's consolidation model \\ of linear poroelasticity}

\author[TUM]{Stephan B. Lunowa\corref{cor}}
\cortext[cor]{Corresponding author}
\ead{stephan.lunowa@tum.de}
\ead[url]{orcid.org/0000-0002-5214-7245}

\author[TUM]{Barbara Wohlmuth}
\ead{wohlmuth@tum.de}
\ead[url]{orcid.org/0000-0001-6908-6015}

\affiliation[TUM]{
    organization={Technical University of Munich, School of Computation, Information and Technology, Department~of~Mathematics},
    addressline={Boltzmannstraße~3},
    postcode={D-85748},
    city={Garching},
    country={Germany}
}

\begin{abstract}
    Biot's consolidation model is a classical model for the evolution of deformable porous media saturated by a fluid and has various interdisciplinary applications.
    While numerical solution methods to solve poroelasticity by typical schemes such as finite differences, finite volumes or finite elements have been intensely studied, lattice Boltzmann methods for poroelasticity have not been developed yet.
    In this work, we propose a novel semi-implicit coupling of lattice Boltzmann methods to solve Biot's consolidation model in two dimensions.
    To this end, we use a single-relaxation-time lattice Boltzmann method for reaction-diffusion equations to solve the Darcy flow and combine it with a recent pseudo-time multi-relaxation-time lattice Boltzmann scheme for quasi-static linear elasticity.
    We employ a multi-grid method for the latter scheme to achieve quasi-optimal computational cost.
    For the coupling between the equations, we develop a centered update scheme, that incorporates both explicit and semi-implicit contributions.
    The numerical results demonstrate that naive (explicit or semi-implicit) coupling schemes lead to instabilities when the poroelastic system is strongly coupled.
    However, the newly developed centered coupling scheme is stable and accurate in all considered cases, even for the Biot--Willis coefficient being one.
    Furthermore, the numerical results for Terzaghi's consolidation problem and a two-dimensional extension thereof highlight that the scheme is even able to capture discontinuous solutions arising from instantaneous loading.
\end{abstract}


\end{frontmatter}

\section{Introduction}
\label{sec:introduction}

The evolution of deformable porous media saturated by a fluid is subject to long-lasting scientific investigations due to its various interdisciplinary applications ranging from geological CO\textsubscript{2} storage \cite{Nordbotten2011} and geothermal energy \cite{Sharmin2023}, over material design for nano-composites \cite{Alkhateb2013} and concrete \cite{Vandamme2013}, to the function of biological tissue \cite{Khaled2003} such as the optic nerve \cite{Causin2014} and heart dynamics \cite{Chapelle2010}.
The theoretical origin dates back to 1923, when Terzaghi developed a description of one-dimensional consolidation \cite{Terzaghi1923}.
In the subsequent decades, the theory was extended by Biot \cite{Biot1941,Biot1941b,Biot1955,Biot1957} leading to the nowadays called Biot's consolidation model.
It describes the two-phase solid-fluid systems based on a linearly elastic, porous solid skeleton that undergoes small deformations, while the fluid flow induced by the deformation of the solid is governed by Darcy's law.
For a general introduction to linear and non-linear poroelasticity, we refer the interested reader to \cite{Detournay1993,Coussy2004}.

Since analytical solutions are barely available, there is a need for numerical methods to solve the coupled system.
Starting from classical finite difference methods \cite{Mercer1999,Gaspar2002,Bean2014} and finite element methods \cite{Lewis1987,Murad1994}, more advanced schemes were developed to improve the approximation of stress and flux and allow for strongly heterogeneous material parameters.
This has led to multiple finite volume methods \cite{Nordbotten2016,Honorio2018,Sokolova2019,Terekhov2020}, and a large class of tailored finite element methods using, e.g., discontinuous Galerkin \cite{Phillips2008,Liu2009}, mixed finite elements \cite{Wheeler2014,Ambartsumyan2020}, or enriched finite elements \cite{Choo2018,Kadeethum2021}.
Furthermore, these methods can account for local mass conservation, avoid locking effects and suppress non-physical pressure oscillations.
For large deformations, the hybrid usage of Lagrangian particles and Eulerian background mesh has been introduced as material point method and extended to poroelasticity \cite{Sulsky1995,Yu2024}.

A substantially different numerical approach are lattice Boltzmann methods (LBMs), which were originally developed for the simulation of fluid dynamics \cite{McNamara1988,Chen1998}.
In contrast to the above mentioned methods, which rely on the approximation of the macroscopic equations, the LBM is based on a mesoscopic description of statistical physics for gas particle dynamics, i.e., the discretization of the Boltzmann equation in physical and velocity space.
While the single-relaxation-time collision operator \cite{Bhatnagar1954} is most popular due to its simplicity, variants based on multiple relaxation times \cite{DHumieres2002} and cumulants \cite{Geier2015} are used for more complex problems to improve accuracy, stability and the number of represented physical parameters.
For a detailed introduction to LBM, we refer to \cite{Krueger2017,Succi2018}.
A major advantage of LBM is the simple algorithmic structure and the straightforward parallelization, especially on GPUs, which make the method well suited for modern, massively parallel and hybrid computing architectures \cite{Schoenherr2011,Ataei2024}.

Nowadays, various LBMs exist for most of the partial differential equations encountered in practice.
In particular, LBMs were quickly adapted for advection-diffusion-reaction equations which have a similar structure as the Navier--Stokes equations.
This was originally recognized for pure diffusion in \cite{WolfGladrow1995}, for advection-diffusion problems in \cite{Flekkoy1993} and for diffusion-reaction problems in \cite{PonceDawson1993}.
Based on these works, the method has been extended in several ways.
Anisotropic diffusion and dispersion was included by a direction-dependent relaxation in \cite{Zhang2002} and using a multiple-relaxation-time (MRT) scheme in \cite{Ginzburg2005,Rasin2005}, while high P\'eclet numbers were treated with an upwind scheme in \cite{Dapelo2021}.
The full advection-diffusion-reaction equation was considered in \cite{Shi2008} using a source-correction scheme involving the temporal derivative, while a simpler moment update was proposed in \cite{Seta2013}.
Furthermore, nonlinear coefficients were incorporated in \cite{Li2015} and combined with a MRT scheme for anisotropy in \cite{Chai2016}.
We refer to \cite{Chai2020} for a detailed analysis of LBM for nonlinear and anisotropic advection-diffusion-reaction equations.
Finally, we note the connection to relaxation systems for general equilibrium and collision operators established in \cite{Simonis2020}, and the fourth-order accuracy achieved in \cite{Geier2017,Chen2024} for a LBM solving diffusion-reaction equations.

Concerning linear elasticity, the first LBM was introduced in \cite{Marconi2003} using two wave equations of the same speed describing the motion of longitudinal and transverse displacements.
This approach was adopted in \cite{OBrien2012} using extended velocity sets in 2D and 3D combined with finite difference schemes and flux limiters to avoid spurious oscillations.
Furthermore, arbitrary Poisson's ratios, and thus different wave speeds, were achieved by an extension to isotropic linear materials in \cite{Murphy2018,Escande2020}.
The realization of several types of boundary conditions has been developed in \cite{Escande2020,Faust2024}.
Another approach is discussed in \cite{Schlueter2018}, where anti-plane shear deformation is assumed to simplify the equations into a scalar wave equation, that is subsequently solved by LBM.
An extension to two wave equations describing dilation and rotation in 2D is presented in \cite{Schlueter2021}.
There, the two equations are solved by two separate LBMs, which are coupled via finite difference schemes.
Recently, a vectorial LBM has been proposed in \cite{Boolakee2025} for linear elastodynamics in two-dimensional domains with periodic or Dirichlet boundary.

While all these methods solve dynamic elasticity, the approximation of quasi-static elasticity by LBM is more intricate.
LBMs are essentially explicit schemes in time, and hence not directly applicable to elliptic equations.
To circumvent this, the addition of a time-dependent damping term was proposed in \cite{Yin2016}, which leads to an iterative LBM that solves the quasi-static equation when approaching the steady-state.
In particular, the method in \cite{Yin2016} combines a vectorial LBM for the components of the displacement with a finite difference approximation of the divergence.
Using the same pseudo-time approach, quasi-static linear elasticity in 2D has been solved in \cite{Boolakee2023} using a MRT LBM on a reduced velocity stencil.
While originally posed as an periodic problem, the treatment of Dirichlet and Neumann boundary conditions is discussed in \cite{Boolakee2023BC}.

To accelerate the convergence to a steady-state solution, multi-grid methods are frequently used for all types of discretization methods.
The fundamental idea of multi-grid methods is that many iterative methods rapidly smooth high frequency error components on a fine grid, while low frequency error components can be represented on a coarser grid, where they become high frequency error components again.
Thus, recursively applying the iterative method on multiple meshes and transferring the corresponding error corrections, this approach achieves a quasi-optimal convergence acceleration.
In the context of LBM, multi-grid methods were first introduced in \cite{Toelke2002,Mavriplis2006} for steady-state solutions of the Navier--Stokes equations.
Multiple extensions for this application were developed.
In \cite{Armstrong2019}, an generalization for MRT LBM was presented, a simplified V-cycle was proposed in \cite{An2023}, and the combination with complex boundaries treated by an immersed-boundary method was discussed in \cite{Gsell2020}.
Moreover, multi-grid LBM for solving second-order elliptic partial differential equations was derived in \cite{Patil2014}.

However, to the best of the authors' knowledge, there have not been developed any LBMs for solving poroelasticity yet.
In this work, we propose a novel coupled LBM to solve Biot's consolidation model in two dimensions.
To this end, we use the classical LBM for reaction-diffusion equations \cite{Seta2013} for the Darcy flow and combine it with the recent pseudo-time LBM for quasi-static linear elasticity \cite{Boolakee2023}.
We employ a multi-grid method for the latter scheme and thus achieve quasi-optimal computational cost.
To resolve the implicit coupling between the equations, we develop a centered update scheme, that incorporates both explicit and semi-implicit contributions.
We demonstrate that naively explicitly or semi-implicitly coupled LBMs lead to instabilities for a strongly coupled problem (Biot--Willis coefficient close to one), while the centered coupling scheme using both contributions is stable and accurate for all Biot--Willis coefficients.
The numerical results for Terzaghi's consolidation problem highlight that the scheme is even able to capture solutions with discontinuities arising from instantaneous loading.

The outline of this paper is as follows.
The governing equations of Biot's consolidation model are presented in \cref{sec:model}.
Next, we introduce the LBMs for flow and elasticity in \cref{sec:LBM}, discussing in detail the different coupling schemes, the multi-grid acceleration and the used initial and boundary conditions.
In \cref{sec:numerics}, we present numerical results for a problem with a manufactured solution, for Terzaghi's consolidation problem as well as for an two-dimensional extension of the latter.
Finally, \cref{sec:conclusion} summarizes the main conclusions.

\section{Governing equations}
\label{sec:model}

In the following two subsections, we first present the system of equations representing Biot's consolidation model in dimensional form and subsequently derive the dimensionless formulation, which will be used in the LBM.

\subsection{Biot's consolidation model}
\label{sec:model:biot}

We consider Biot's model of consolidation for a linearly elastic and porous solid, which is saturated by a slightly compressible viscous fluid \cite{Biot1941b,Biot1957,Detournay1993,Coussy2004}.
For given final time $t_f \in (0, \infty)$ and a simply connected, Lipschitz-bounded domain $\mathcal{D} \in \mathbb{R}^d$, we define the time-space cylinder $\mathcal{Q} = (0, t_f) \times \mathcal{D}$.
Note that we discuss in this section the general model with dimension $d \in \Set{1, 2, 3}$, whereas we restrict ourselves to $d = 2$ thereafter for the LBM.
Furthermore, the displacement and stress of the solid are denoted by $\vec{\eta} : \mathcal{Q} \to \mathbb{R}^d$ and $\tensor{\upsigma} : \mathcal{Q} \to \mathbb{R}^{d \times d}$, respectively, while the fluid pressure and Darcy velocity are denoted by $p : \mathcal{Q} \to \mathbb{R}$ and $\vec{u} : \mathcal{Q} \to \mathbb{R}^d$, respectively.
Assuming an isotropic material, the system of partial differential equations reads
\begin{align}
    - \nabla\cdot\tensor{\upsigma} + \alpha \nabla p &= \vec{f} &&\text{in }\mathcal{Q}, \label{eq:elasticity}\\
    \tensor{\upsigma} &= \lambda \big(\nabla\cdot\vec{\eta}\big)\tensor{I} + \mu \big(\nabla\vec{\eta} + \nabla\vec{\eta}^\top\big) &&\text{in }\mathcal{Q}, \label{eq:stress}\\
    \partial_{t} \big(c_0 p + \alpha \nabla\cdot\vec{\eta}\big) + \nabla\cdot\vec{u} &= s &&\text{in }\mathcal{Q}, \label{eq:flow}\\
    \vec{u} &= - \kappa \big(\nabla p - \vec{g}\big) &&\text{in }\mathcal{Q}. \label{eq:darcy}
\end{align}
Here, the first two equations describe the quasi-static force balance in the deformed solid, while the last two equations describe the diffusive Darcy flow of the fluid due to the deformation.
The so-called Biot--Willis parameter $\alpha \in (0, 1]$ represents the interaction strength between solid and fluid ranging from no coupling in the limit $\alpha \to 0$ to strong coupling for $\alpha$ close to 1.
Moreover, $\lambda$ and $\mu$ denote the first and second Lam\'e coefficients of the solid, $c_0$ is the specific storage coefficient, and $\kappa$ represents the viscosity-reduced permeability.
Finally, the right-hand side $\vec{f}$ and $\vec{g}$ are given volume forces for solid and fluid, respectively, while $s$ is a volumetric fluid source or sink.

The system has to be closed by appropriate initial and boundary conditions for the primal unknowns $\vec{\eta}$ and $p$.
We omit the details here and will specify the conditions used for the numerical examples in \cref{sec:numerics} and their LBM implementations in \cref{sec:LBM:BCs}.
Note that a derivation of this model by homogenization can be found in \cite{Auriault1977}, while a rigorous mathematical analysis of this problem is presented in \cite{Showalter2000}.

\subsection{Dimensionless system of equations}
\label{sec:model:dimensionless}

For the dimensionless formulation, we introduce a characteristic length $L$, time $T$, scalar displacement $H$ and pressure $P = c_0^{-1}$.
With this reference system, the dimensionless quantities are defined by
\begin{equation*}
    \hat{\vec{x}} = \frac{\vec{x}}{L}, \qquad
    \hat{t} = \frac{t}{T}, \qquad
    \hat{\vec{\eta}} = \frac{\vec{\eta}}{H}, \qquad
    \hat{\tensor{\upsigma}} = \frac{c_0 L\tensor{\upsigma}}{H}, \qquad
    \hat{p} = c_0 p, \qquad
    \hat{\vec{u}} = \frac{T \vec{u}}{L},
\end{equation*}
while the dimensionless derivatives satisfy $\hat{\nabla} = L \nabla$ and $\partial_{\hat{t}} = T \partial_{t}$.
Accordingly, the scaled time-space cylinder is given by $\hat{\mathcal{Q}} = \Set{ (\hat{t}, \hat{\vec{x}}) \given (t, \vec{x}) \in \mathcal{Q} }$.
Inserting the dimensionless quantities into the Biot system \crefrange{eq:elasticity}{eq:darcy} yields
\begin{align}
    - \hat\nabla\cdot \hat{\tensor{\upsigma}} + \varepsilon\alpha \hat\nabla \hat{p} &= \varepsilon \hat{\vec{f}\,}\! &&\text{in }\hat{\mathcal{Q}}, \label{eq:elasticity:dim}\\
    \hat{\tensor{\upsigma}} &= \hat\lambda \big(\hat\nabla\cdot \hat{\vec{\eta}}\big)\tensor{I} + \hat\mu \big( \hat\nabla \hat{\vec{\eta}} + \hat\nabla \hat{\vec{\eta}}^\top\big) &&\text{in }\hat{\mathcal{Q}}, \label{eq:stress:dim}\\
    \partial_{\hat{t}} \left(\hat{p} + \frac{\alpha}{\varepsilon} \hat\nabla\cdot \hat{\vec{\eta}}\right) + \hat\nabla\cdot \hat{\vec{u}} &= \hat{s} &&\text{in }\hat{\mathcal{Q}}, \label{eq:flow:dim}\\
    \hat{\vec{u}} &= - \hat\kappa \big(\hat\nabla \hat{p} - \hat{\vec{g}}\big) &&\text{in }\hat{\mathcal{Q}}, \label{eq:darcy:dim}
\end{align}
where the dimensionless parameters are
\begin{equation*}
    \varepsilon = \frac{L}{H}, \qquad
    \hat{\lambda} = c_0 \lambda, \qquad
    \hat{\mu} = c_0 \mu, \qquad
    \hat{\kappa} = \frac{T \kappa}{c_0 L^2},
\end{equation*}
and the dimensionless forces and source are given by
\begin{equation*}
    \hat{\vec{f}\,}\! = c_0 L \vec{f}, \qquad
    \hat{\vec{g}} = c_0 L \vec{g}, \qquad
    \hat{s} = T s.
\end{equation*}
In the following, we drop the hats for readability.

\paragraph{Remark on the incompressible limit cases}
For an incompressible fluid ($c_0 = 0$), Biot's consolidation model is governed by a degenerate parabolic system as the temporal derivative of the pressure vanishes.
Then, the above scaling is impossible, but one can use an appropriate scaling, e.g. by replacing $c_0$ by $\mu$ for the definition of the dimensionless quantities, leading to well-defined equations.
Since the pressure equation is elliptic in that case, standard LBM cannot be used directly, but one also needs a pseudo-timestepping for the pressure equation, similar to the one for linear elasticity presented below.
For simplicity, we restrict ourselves to the compressible case ($c_0 > 0$) in this article.

For an incompressible solid ($\lambda \to \infty$), the elastic sub-system must be augmented by the additional constraint $\nabla\cdot\vec{\eta} = 0$.
However, the LBM for linear elasticity presented below becomes unstable in this limit \cite{Boolakee2023}.
In the context of finite element and volume methods, remedies such as the total pressure formulation are known and might be applied analogously to the presented LBM.
However, this goes beyond the scope of this article, so that we only consider compressible solids from here on.

\section{The lattice Boltzmann method for Biot's consolidation model}
\label{sec:LBM}

In this section, we propose the LBM for Biot's consolidation model by coupling LBM for the diffusive fluid flow \cite{Seta2013} with LBM for two-dimensional linear elasticity \cite{Boolakee2023,Boolakee2023BC}.
To this end, we introduce notation conventions and basic definitions first and subsequently present the separate LBMs for diffusive flow and elasticity.
Then, we discuss the coupling of the two methods, explain the multi-grid method accelerating the pseudo-timestepping, and finally elaborate on the LBM implementation of the used initial and boundary conditions.

\subsection{Notation and definitions for lattice Boltzmann methods}
\label{sec:LBM:basics}

\begin{figure}[tb]
    \centering
    \includegraphics{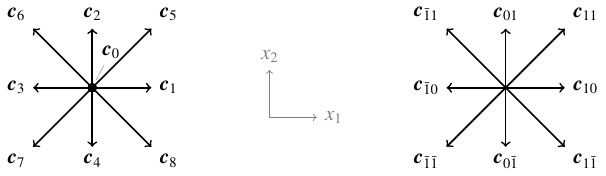}
    \caption{D2Q9 velocity set $\vec{c}_i$ with linear indices (left) and D2Q8 velocity set $\vec{c}_{ij}$ with 2D Miller indices, i.e., $\bar 1 = -1$ (right).}\label{fig:velocities}
\end{figure}

In the following, we summarize all basics of LBM necessary for this work, for a detailed introduction we refer to \cite{Krueger2017,Succi2018}.
The fundamental principle of LBM is the discretization of the Boltzmann equation in physical and velocity space.
For this purpose, the main variable is the distribution function $f$, which denotes the probability of material particles at time $t$ and position $\vec{x}$ to have the velocity $\vec{c}$.
The velocity space is discretized using a finite set of velocities $\vec{c}_i$ and corresponding weights $w_i$, such that $f_i(t, \vec{x}) = f(t, \vec{x}, \vec{c}_i)$.
These sets are typically derived based on the required mass, momentum and energy conservation using a Hermite series expansion.
Typically, the velocity sets are denoted by D$d$Q$q$, where $d$ is the spatial dimension and $q$ the number of discrete velocities.
For the two-dimensional LBMs presented in the following, we utilize the classical D2Q9 velocity set as well as a reduced D2Q8 velocity set without rest velocity $\vec{c}_0 = \vec{0}$, see \cref{fig:velocities,tab:velocities}.

\begin{table}[tb]
    \renewcommand\arraystretch{1.25}
    \newcommand*\Mvec[2]{\begin{psmallmatrix*}[r]#1\\#2\end{psmallmatrix*}}
    \centering
    \caption{Discrete velocities $\vec{c}_i$ and weights $w_i$ with linear and 2D Miller indices, i.e., $\bar 1 = -1$.}\label{tab:velocities}
    \vspace{1ex}
    \begin{tabular}{l|*{9}{C}}
    \toprule
    Linear index &           0 &           1 &           2 &            3 &            4 &          5  &            6 &             7 &            8 \\
    \midrule
    Miller index &          00 &          10 &          01 &      \bar 10 &      0\bar 1 &          11 &      \bar 11 &  \bar 1\bar 1 &      1\bar 1 \\
    \midrule
    Velocity     & \Mvec{0}{0} & \Mvec{1}{0} & \Mvec{0}{1} & \Mvec{-1}{0} & \Mvec{0}{-1} & \Mvec{1}{1} & \Mvec{-1}{1} & \Mvec{-1}{-1} & \Mvec{1}{-1} \\
    \midrule
    Weight       &         4/9 &         1/9 &         1/9 &          1/9 &          1/9 &        1/36 &         1/36 &          1/36 &         1/36 \\
    \bottomrule
    \end{tabular}
\end{table}

Additionally, the two-dimensional spatial domain $\mathcal{D}$ is discretized by a uniform lattice of width $\Delta x$, while the time interval $(0, t_f)$ is discretized by a uniform time step $\Delta t = t_f / N_T$ with $N_T \in \N$.
This leads to a discrete lattice speed of sound $c_s = \Delta x / (\sqrt{3} \Delta t)$.
Altogether, the general lattice Boltzmann equation then reads
\begin{equation}
    f_i(t+\Delta t, \vec{x}+\Delta t\vec{c}_i) - f_i(t, \vec{x}) = \Delta t \, \Omega_i(t, \vec{x}), \label{eq:LBM}
\end{equation}
where the left-hand side represents an advection along the discrete velocity $\vec{c}_i$, the so-called streaming step, whereas the right-hand side term, called collision operator, contains all interactions such as relaxation to equilibrium due to particle collision as well as external forces and sources.
Note that the lattice Boltzmann equation can be split in two parts: the collision step
\begin{equation}
    f_i^*(t, \vec{x}) = f_i(t, \vec{x}) + \Delta t \, \Omega_i(t, \vec{x}), \label{eq:collision}
\end{equation}
which is local, and the streaming step
\begin{equation}
    f_i(t+\Delta t, \vec{x}+\Delta t\vec{c}_i) = f_i^*(t, \vec{x}), \label{eq:streaming}
\end{equation}
which moves each distribution function one lattice node further in their velocity direction.
This is the fundamental reason for the algorithmic simplicity and for the straightforward and well-scaling parallelization.

Even in the absence of external interactions, the collision operator is very complicated as it involves all types of microscopic particle interactions.
However, the most simple approximation, which ensures mass and momentum conservation, is given by the Bhatnagar--Gross--Krook collision operator \cite{Bhatnagar1954}
\begin{equation}
    \Omega_i = -\omega \big(f_i - f_i^\text{eq}\big), \label{eq:SRT}
\end{equation}
where $\omega > 0$ is the relaxation rate towards the equilibrium state $f_i^\text{eq}$ that depends on the moments of all distribution functions $f_j$ according to the equation of interest,
e.g., on density $\rho = \sum_i f_i$ and momentum $\rho\vec{u} = \sum_i f_i \vec{c}_i$ in case of the Navier--Stokes equations.
Note that here and in the following, the arguments $(t,\vec{x})$ are omitted whenever the meaning is unambiguous.
The Bhatnagar--Gross--Krook collision operator will be used below for the diffusive fluid flow.
In this application, the relaxation rate is directly related to the diffusion coefficient.

In contrast to this single-relaxation-time scheme, multiple-relaxation-time (MRT) collision operators allow for a larger number of free parameters, increased stability and higher accuracy.
This is necessary to formulate a LBM for linear elasticity.
For MRT schemes, the distribution functions are transformed into linear independent moments, which are individually relaxed at different rates and then transformed back to the distributions \cite{DHumieres2002}.
Denoting the vector of the distribution functions for all discrete velocities by $\vec{f} = (f_0, \dots, f_8)^\top$, the transformation into moments can be described by multiplication by an invertible matrix $\tensor{M}$.
Analogously, we represent the individual relaxation rates as a diagonal matrix $\tensor{\omega}$ and the equilibrium moment vector by $\vec{m}^\text{eq}$.
The vector of MRT collision operators for all discrete velocities is then given by
\begin{equation}
    \vec{\Omega} = -\tensor{M^{-1}} \tensor{\omega} \big(\tensor{M} \vec{f} - \vec{m}^\text{eq}\big) . \label{eq:MRT}
\end{equation}
The specific moments and relaxation rates depend again on the parameters of the considered equations.
More details will be given below during the discussion of LBM for linear elasticity.
Finally, note that the Bhatnagar--Gross--Krook collision operator is a special case of the MRT one when $\tensor{\omega} = \omega \tensor{I}$ is chosen.

\subsection{The lattice Boltzmann method for diffusive fluid flow}
\label{sec:LBM:flow}

In the following, we present the required steps for the LBM solving the diffusive fluid flow \cref{eq:flow:dim,eq:darcy:dim}.
Combining these equations and reordering the terms, we obtain the diffusion-reaction equation
\begin{align}
    \partial_{t} p - \kappa \Delta p &= s_\text{eff} = s - \kappa \nabla\cdot\vec{g} - \frac{\alpha}{\varepsilon} \partial_{t} \nabla\cdot \vec{\eta} &&\text{in }\mathcal{Q}. \label{eq:flow:combined}
\end{align}
Subsequently, the LBM for the diffusion-reaction equation proposed in \cite{Seta2013} is discussed, while the detailed computation of the effective source term $s_\text{eff}$ is deferred to \cref{sec:LBM:coupling}.
Following the presentation in \cite[Chpt.~8]{Krueger2017}, this LBM is based on the single-relaxation-time scheme \cref{eq:LBM,eq:SRT} and uses the D2Q9 velocity set.
Denoting the distribution functions for the pressure by $f_i$ with linear index $i \in \Set{0, 1, \dots, 8}$,
the pressure $p$ is defined by the zeroth moment of the distribution functions corrected for the effective source $s_\text{eff}$ as
\begin{equation}
    p = \sum_{i=0}^8 f_i + \frac{\Delta t}{2} s_\text{eff} .\label{eq:pressure}
\end{equation}
Since \cref{eq:flow:combined} is purely diffusive, the equilibrium distribution functions are defined by
\begin{equation}
    f^\text{eq}_i = w_i p .\label{eq:pressure:equilibrium}
\end{equation}
The single-relaxation-time collision step combining \cref{eq:collision,eq:SRT} then reads
\begin{equation}
    f_i^* = f_i - \omega \Delta t \big( f_i - f^\text{eq}_i \big) + \Delta t S_i ,\label{eq:pressure:collision}
\end{equation}
where the contribution by the effective source is given by
\begin{equation}
    S_i = \left(1-\frac{\omega}{2}\right) w_i s_\text{eff} .\label{eq:pressure:collision:source}
\end{equation}
Furthermore, the relaxation rate $\omega$ must satisfy $\omega \in (0, 2)$ for stability and is related to the permeability via
\begin{equation}
    \kappa = c_s^2 \left(\frac{1}{\omega} - \frac{\Delta t}{2}\right) .\label{eq:pressure:rate}
\end{equation}
Note that this is the standard relation between relaxation rate $\omega$ and diffusion coefficient $\kappa$ in LBM.
In the end of each time step, the streaming step of \cref{eq:streaming} is performed.
Overall, this LBM yields a pressure solution that is second-order consistent in time and space \cite{Seta2013}.

\subsection{The lattice Boltzmann method for quasi-static linear elasticity}
\label{sec:LBM:elasticity}

Next, we present the necessary steps for the LBM solving the quasi-static linear elasticity \cref{eq:elasticity:dim,eq:stress:dim}.
The method has been developed in \cite{Boolakee2023,Boolakee2023BC} and uses a MRT scheme for the D2Q8 velocity set.
Here, we adapt the description to the context of Biot's consolidation model, where an additional parameter-like time-dependence emerges.
In particular, the method solves \cref{eq:elasticity:dim,eq:stress:dim} using an pseudo-timestepping method, where \cref{eq:elasticity:dim} is replaced by
\begin{align}
    \frac{1}{\varepsilon} \partial_{\tau}\vec{\eta} - \nabla\cdot \tensor{\upsigma} &= \varepsilon \vec{f}_\text{eff} = \varepsilon \big(\vec{f} - \alpha \nabla p\big) &&\text{for }(\tau, t, \vec{x}) \in (0, \infty) \times \mathcal{Q}, \label{eq:elasticity:pseudo}
\end{align}
which is equivalent when the stationary solution is reached (for $\tau \to \infty$).
Note that the division by $\varepsilon$ in front of the pseudo-time derivative is a result of the viscous scaling leading to a pseudo-time step $\Delta\tau = 1$.
As before, we defer the computation of the effective force $\vec{f}_\text{eff}$ to \cref{sec:LBM:coupling}, and elaborate on the general method performed within each time step for multiple pseudo-time steps.

The moments for the MRT scheme are linear combinations of the raw countable moments \cite{Geier2015}
\begin{equation}
    m_{ab} = \sum_{i=-1}^1 \sum_{j=-1}^1 i^a j^b g_{ij} ,\label{eq:moments}
\end{equation}
using the two-dimensional Miller index $ab \in \Set{ 10, 01, 11, 20, 02, 12, 21, 22 }$, cf. also \cref{tab:moments:inverse} for the inverse transformation.
Note that the number of independent moments coincides with the size of the velocity set, i.e., there is no independent zeroth moment $m_{00} = \sum_{i,j} g_{ij}$ due to the absence of the rest velocity $\vec{c}_{00}$.
Besides the raw moments $m_{ab}$ with $ab \in \Set{ 10, 01, 11, 12, 21 }$, we use the combined second-order moments
\begin{equation}
    m_s = m_{20} + m_{02} ,\qquad
    m_d = m_{20} - m_{02} ,\label{eq:moments:stress}
\end{equation}
which are directly related to the volumetric and deviatoric components of the stress.
Additionally, the raw fourth-order moment $m_{22}$ is replaced by
\begin{equation}
    m_f = m_{22} + \frac{1}{12(\lambda + \mu) - 4} m_s ,\label{eq:moments:BC}
\end{equation}
which is required for imposing Neumann boundary conditions \cite{Boolakee2023BC}, see \cref{sec:LBM:BCs}.

Interpreting the LBM as a Strang splitting scheme, the effective force is split into two sub-steps of the collision step \cite{Dellar2013}, and then applied as follows
\begin{equation}
    \begin{pmatrix} \bar{m}_{10} \\ \bar{m}_{01} \end{pmatrix} = \begin{pmatrix} m_{10} \\ m_{01} \end{pmatrix} + \frac{\varepsilon^2}{2} \vec{f}_\text{eff} ,\qquad
    \begin{pmatrix} m^*_{10} \\ m^*_{01} \end{pmatrix} = \begin{pmatrix} \bar{m}_{10} \\ \bar{m}_{01} \end{pmatrix} + \frac{\varepsilon^2}{2} \vec{f}_\text{eff} .\label{eq:moments:halfForce}
\end{equation}
Note that this reflects a conservation of the first moments, since no relaxation is applied ($\omega_{10} = \omega_{01} = 0$).
All other moments $m_\nu$, $\nu \in \Set{ s, d, 11, 12, 21, f }$, are relaxed according to the MRT scheme \cref{eq:MRT}, viz.,
\begin{equation}
    m^*_\nu = m_\nu - \omega_\nu (m_\nu - m^\text{eq}_\nu), \label{eq:moments:collision}
\end{equation}
where the relaxation rates $\omega_\nu$ and equilibrium moments $m^\text{eq}_\nu$ are given in \cref{tab:moments}.
In particular, note that the results after half the collision step are denoted by $\bar{m}_\nu = (m_\nu + m^*_\nu) / 2$.

\begin{table}[tb]
    \renewcommand\arraystretch{1.25}
    \newcommand*\WS{2 / \big(3(\lambda+\mu) + 1\big)}
    \newcommand*\WD{2 / (6\mu + 1)}
    \centering
    \caption{Relaxation rates $\omega_\nu$ and equilibrium moments $m^\text{eq}_\nu$ for quasi-static linear elasticity using 2D Miller indices.}
    \label{tab:moments}
    \vspace{1ex}
    \begin{tabular}{l|*{8}{C}}
        \toprule
        Index $\nu$                                     & 10 & 01 &   s &   d &  11 &                      12 &                      21 & f \\
        \midrule
        Relaxation rate $\omega_\nu$                  &  0 &  0 & \WS & \WD & \WD &                       1 &                       1 & 1 \\
        \midrule
        Equilibrium moment $m^\text{eq}_\nu$ &  - &  - &   0 &   0 &   0 & \bar{m}_{10}/3 & \bar{m}_{01}/3 & 0 \\
        \bottomrule
    \end{tabular}
\end{table}

Furthermore, the relaxation rates $\omega_d$ and $\omega_{11}$ are the same in order to obtain isotropic behavior of the deviatoric stress.
Indeed, it has been shown in \cite{Boolakee2023BC} by asymptotic expansion that the displacement $\vec{\eta}$ and the stress $\tensor{\upsigma}$ are given in the limit $\tau \to \infty$ by
\begin{equation}
    \vec{\eta} = \begin{pmatrix} \bar{m}_{10} \\ \bar{m}_{01} \end{pmatrix} ,\qquad
    \tensor{\upsigma} = -\frac{1}{2} \begin{pmatrix} \bar{m}_s + \bar{m}_d & 2\bar{m}_{11} \\ 2\bar{m}_{11} & \bar{m}_s - \bar{m}_d \end{pmatrix} ,\label{eq:displacement}
\end{equation}
respectively, where the involved moments after half the collision can be computed as
\begin{equation}
    \bar{m}_s = \frac{3(\lambda + \mu)}{3(\lambda + \mu) + 1} m_s ,\qquad
    \bar{m}_d = \frac{6\mu}{6\mu + 1} m_d ,\qquad
    \bar{m}_{11} = \frac{6\mu}{6\mu + 1} m_{11} .\label{eq:moments:bared}
\end{equation}

The post-collision distribution functions are obtained from the back-transformation of the post-collision moments by inverting \cref{eq:moments,eq:moments:stress,eq:moments:BC}.
First, note that this leads to
\begin{equation}
    m^*_{22} = - \frac{1}{12(\lambda+\mu) - 4} m^*_s = - \frac{1}{12(\lambda+\mu) + 4} m_{s},
    \label{eq:moments:collision:22}
\end{equation}
which is non-singular even for $\lambda + \mu = 1/3$ (alternatively, in this case, the relaxation rate $\omega_f$ can be changed without influencing the leading-order behavior \cite{Boolakee2023BC}).
Then, the post-collision distribution functions are given by
\begin{equation}
    g^*_{ij} = \sum_\nu c_{ij}^\nu m^*_\nu ,\label{eq:moments:inverse}
\end{equation}
where the coefficients $c_{ij}^\nu$ are stated in \cref{tab:moments:inverse}.
In the end of each pseudo-time step, the streaming step of \cref{eq:streaming} is performed for the distribution functions $g_{ij}$.
Overall, this LBM yields displacement and stress solutions which are second-order consistent in space for pseudo time $\tau \to \infty$ (in every time step) \cite{Boolakee2023}.

\begin{table}[tbp]
    \renewcommand\arraystretch{1.5}
    \newcommand*{\h}{\frac12}
    \newcommand*{\q}{\frac14}
    \centering
    \caption{Coefficients $c_{ij}^\nu$ of the back-transformation from moments $m^*_\nu$ to distribution functions $g^*_{ij}$ using 2D Miller indices, i.e., $\bar 1 = -1$.}
    \label{tab:moments:inverse}
    \vspace{1ex}
    \begin{tabular}{C|*{8}{R}}
        \toprule
        {}_{\textstyle ij} \diagdown {}^{\textstyle\nu} & 01 & 10 & 11 & s & d & 12 & 21 & 22 \\
        \midrule
                  10 &  \h &   0 &   0 & \q &  \q & -\h &   0 & -\h \\
                  01 &   0 &  \h &   0 & \q & -\q &   0 & -\h & -\h \\
             \bar 10 & -\h &   0 &   0 & \q &  \q &  \h &   0 & -\h \\
             0\bar 1 &   0 & -\h &   0 & \q & -\q &   0 &  \h & -\h \\
                  11 &   0 &   0 &  \q &  0 &   0 &  \q &  \q &  \q \\
             \bar 11 &   0 &   0 & -\q &  0 &   0 & -\q &  \q &  \q \\
        \bar 1\bar 1 &   0 &   0 &  \q &  0 &   0 & -\q & -\q &  \q \\
             1\bar 1 &   0 &   0 & -\q &  0 &   0 &  \q & -\q &  \q \\
        \bottomrule
    \end{tabular}
\end{table}

\subsection{The coupling of the lattice Boltzmann methods}
\label{sec:LBM:coupling}

To obtain the LBM for Biot's consolidation model, we couple the above discussed LBMs for diffusive flow and for quasi-static linear elasticity.
To this end, in each time step we solve the quasi-static elasticity problem using $N_E \in \N$ pseudo-time steps, followed by
one regular update of the LBM for diffusive flow.

The effective force term $\vec{f}_\text{eff} = \vec{f} - \alpha \nabla p$ in the elasticity \cref{eq:elasticity:pseudo} involves the gradient of the pressure.
Following \cite{Lee2005,Lee2006}, this gradient can be approximated up to second-order errors in space by the central difference quotient
\begin{equation}
    \nabla^c p(\vec{x})
    = \sum_{i=0}^8 \frac{w_i\vec{c}_i \big(p(\vec{x}+\Delta t\vec{c}_i) - p(\vec{x}-\Delta t\vec{c}_i)\big)}{2c_s^2\Delta t}
    = \sum_{i=0}^8 \frac{w_i\vec{c}_i p(\vec{x}+\Delta t\vec{c}_i)}{c_s^2\Delta t} .\label{eq:pressure:gradient:alternative}
\end{equation}
Note that the second equality only holds for symmetric velocity sets, such as the D2Q9 used here.
Alternatively, one can use the first moments of the pressure distribution functions $f_i$ to obtain the pressure gradient \cite[Chpt.~8]{Krueger2017}.
Accounting for the contribution of the source term, this yields
\begin{equation}
    \nabla^* p = -\frac{\omega}{c_s^2} \sum_{i=0}^8 f_i \vec{c}_i + \frac{1}{2}\nabla^c s_\text{eff} .\label{eq:pressure:gradient}
\end{equation}
For the numerical examples below, we observe the same order of convergence for both methods, but considerably smaller errors when using \cref{eq:pressure:gradient}, so that we only present results of the latter in the following.

Irrespective of the used approximation, the gradient depends on the effective source term $s_\text{eff} = s - \kappa \nabla\cdot\vec{g} - \frac{\alpha}{\varepsilon} \partial_{t} \nabla\cdot\vec{\eta}$ in the flow \cref{eq:flow:combined}.
While the first and second right-hand side terms are given explicitly, the last one does depend again on the solution of the elastic equation.
To eliminate this cyclic self-dependence, one can use a fully explicit scheme
\begin{equation}
    s_\text{eff}^\text{ex}(t) = s(t) - \kappa \nabla\cdot\vec{g}(t) - \frac{\alpha}{\varepsilon} \partial_{t}^- \nabla\cdot\vec{\eta}(t-\Delta t) \label{eq:source:explicit}
\end{equation}
where $\partial_{t}^- f(t) = \big(f(t) - f(t-\Delta t)\big) / \Delta t$ denotes the temporal backward difference, while the divergence of the displacement is given according to \cref{eq:stress,eq:displacement} by
\begin{equation}
    \nabla\cdot\vec{\eta} = - \frac{1}{2(\lambda+\mu)} \bar{m}_s .\label{eq:LBM:divergence}
\end{equation}
However, this scheme is unstable for $\alpha$ close to 1, as shown numerically in \cref{sec:numerics:manufactured}.
Instead, we propose a semi-implicit scheme that utilities the elasticity solution of the current physical time $t$, but at the previous pseudo-time step $\tau-\Delta\tau$ (which is iterated), viz.
\begin{equation}
    s_\text{eff}^\text{im}(\tau, t) = s(t) - \kappa \nabla\cdot\vec{g}(t) - \frac{\alpha}{\varepsilon} \Big( (1 - r)\, \partial_{t}^- \nabla\cdot\vec{\eta}(t-\Delta t) + r\, \partial_{t}^- \nabla\cdot\vec{\eta}(\tau-\Delta\tau, t) \Big)  ,\label{eq:source:implicit}
\end{equation}
where $r \in [0, 1]$ is the ratio between explicit and semi-implicit approximation.
Note that only the rightmost term needs to be recomputed in each pseudo-time step, while the others are computed at the beginning of each time step.
Therefore, each time step is initialized by setting $g_{ij}(-1,t) = g_{ij}(0,t) = g_{ij}(N_E, t-\Delta t)$,
and after the $N_E$ steps one sets $g_{ij}(t) = g_{ij}(N_E,t)$ and $\nabla\cdot\vec{\eta}(t) = \nabla\cdot\vec{\eta}(N_E, t)$.
Altogether, one obtains \cref{alg:LBM} based on the \cref{alg:LBM:elasticity,alg:LBM:flow} for elasticity and flow updates.
While the discussion above only considers the interior of the domain, the algorithms already include the application of initial and boundary conditions, which will be discussed in detail in \cref{sec:LBM:BCs}.

\begin{algorithm}[tbp]
    \caption{Semi-implicit LBM for poroelasticity.}
    \label{alg:LBM}
    \begin{algorithmic}[1]
    \State Set $t = 0$ and initialize the distribution functions $f_i$ and $g_{ij}$                     \Comment{\cref{eq:ICs,eq:ICs:moments}}
    \For{$i = 1, \dots, N_T$}                                                                                  \label{alg:for}
        \LComment{Elasticity update}
        \State Set $\tau = 0$ and compute the explicit part of $s_\text{eff}^\text{im}(t)$             \Comment{\cref{eq:source:implicit}}
        \For{$j = 1, \dots, N_E$}
            \State Compute $s_\text{eff}^\text{im}(\tau,t)$ and $p(\tau,t)$ \Comment{\cref{eq:pressure,eq:source:implicit}}
            \State Compute $\vec{f}_\text{eff}(\tau,t)$ from $\nabla p(\tau,t)$              \Comment{\cref{eq:pressure:gradient}}
            \State Compute $g_{ij}(\tau+\Delta\tau,t) = \Call{ELASTIC\_STEP}{g_{ij}(\tau,t), \vec{f}_\text{eff}(\tau,t)}$ \Comment{\cref{alg:LBM:elasticity}}
            \State Compute $\nabla\cdot\vec{\eta}(\tau+\Delta\tau,t)$                          \Comment{\cref{eq:moments:bared,eq:LBM:divergence}}
            \State Advance pseudo time $\tau$ by $\Delta\tau$
        \EndFor
        \State Compute $\vec{\eta}(t)$ and $\tensor{\upsigma}(t)$                                                      \Comment{\cref{eq:moments:bared,eq:displacement}}
        \LComment{Flow update}
        \State Compute $s_\text{eff}^\text{im}(t)$ and $p(\tau,t)$            \Comment{\cref{eq:pressure,eq:source:implicit}}
        \State Compute $f_i(t+\Delta t) = \Call{FLOW\_STEP}{f_i(t), s_\text{eff}^\text{im}(t)}$ \Comment{\cref{alg:LBM:flow}}
        \State Advance time $t$ by $\Delta t$
    \EndFor
    \end{algorithmic}
\end{algorithm}

\begin{algorithm}[tbp]
    \caption{LBM time step for quasi-static linear elasticity.}
    \label{alg:LBM:elasticity}
    \begin{algorithmic}[1]
        \Function{ELASTIC\_STEP}{$g_{ij}$, $\vec{f}_\text{eff}$}
            \State Compute the moments $m_\nu$ from $g_{ij}$                                         \Comment{\cref{eq:moments,eq:moments:stress,eq:moments:BC}}
            \State Apply half-forcing (first half of Strang splitting)                                     \Comment{\cref{eq:moments:halfForce}}
            \State Compute the equilibrium moments $m^\text{eq}_\nu$                                     \Comment{\cref{tab:moments}}
            \State Apply half-forcing (second half of Strang splitting)                                    \Comment{\cref{eq:moments:halfForce}}
            \State Compute the post-collision moments $m^*_\nu$                                \Comment{\cref{eq:moments:collision,tab:moments}}
            \State Back-transform to distributions $g^*_{ij}$ from $m^*_\nu$ \Comment{\cref{eq:moments:collision:22,eq:moments:inverse}}
            \State Pseudo-time streaming of the distributions: $g_{ij} = \Call{STREAM}{g^*_{ij}}$ \Comment{\cref{eq:streaming}}
            \State Apply the (nontrivial) boundary conditions for linear elasticity                                 \Comment{Implementation details in \cref{sec:LBM:BCs}}
            \State \Return $g_{ij}$
        \EndFunction
    \end{algorithmic}
\end{algorithm}

\begin{algorithm}[tbp]
    \caption{LBM time step for diffusive fluid flow.}
    \label{alg:LBM:flow}
    \begin{algorithmic}[1]
        \Function{FLOW\_STEP}{$f_i$, $s_\text{eff}$}
            \State Compute the equilibrium distributions $f^\text{eq}_i$                        \Comment{\cref{eq:pressure,eq:pressure:equilibrium}}
            \State Compute the post-collision distributions $f^*_i$                   \Comment{\cref{eq:pressure:collision,eq:pressure:collision:source}}
            \State Streaming of the distributions: $f_i = \Call{STREAM}{f^*_i}$ \Comment{\cref{eq:streaming}}
            \State Apply the (nontrivial) boundary conditions for the pressure                 \Comment{Implementation details in \cref{sec:LBM:BCs}}
            \State \Return $f_i$
        \EndFunction
    \end{algorithmic}
\end{algorithm}

Note that an adaptive version of this algorithm using a variable number of pseudo-time steps $N_E$ might be preferable for solutions approaching a stationary state.
From an algorithmic point of view, such a change is trivial, as one simply replaces \cref{alg:for} by a while loop using the stopping criterion.
However from the mathematical theory, it remains unclear, what would be an appropriate stopping criterion.
Without known convergence order for the elastic sub-system, a simple criterion using the difference between consecutive iterations could easily fail.
Therefore, we restrict ourselves here to the study of the algorithm with fixed number of pseudo-time steps.

However, the underlying elasticity problem is elliptic, and thus can be efficiently treated by multi-grid methods, see e.g. \cite{Mavriplis2006,Patil2014}.
In the following section, we propose such a method for the elasticity LBM discussed above.
Thereby, the number of pseudo-time steps becomes essentially independent of the spatial discretization yielding an quasi-optimal algorithm.

Furthermore, with respect to the dimensionless formulation discussed in \cref{sec:model:dimensionless},
note that in the LBM setting one typically chooses the characteristic scales $T \sim \Delta t_\text{phys}$ and $L \sim \Delta x_\text{phys}$ to obtain step sizes $\Delta t_\text{LBM} = \Delta x_\text{LBM} = 1$.
Together with a fixed $H \sim 1$, this also yields a small length ratio $\varepsilon \sim \Delta x_\text{phys}$ required for the convergence of the LBM for elasticity \cite{Boolakee2023}.
Finally, using the viscous scaling $\Delta t_\text{phys} \sim (\Delta x_\text{phys})^2$ yields constant material parameters, but decreasing forces and source.

\subsection{Acceleration by a multi-grid method}
\label{sec:LBM:MG}

To accelerate the convergence of the LBM for quasi-static linear elasticity, we use a multi-grid method, as discussed in the following.
Note that when the steady-state solution is reached, applying one LBM step results in the same populations.
Thus, the residual of the LBM for quasi-static linear elasticity is given by
\begin{equation}
    r_{ij} = g_{ij} - g^{**}_{ij} , \label{eq:LBM:residual}
\end{equation}
where $g_{ij}^{**}$ is obtained by applying one pseudo-time step (\cref{alg:LBM:elasticity}), viz.
\begin{equation}
    g_{ij}^{**}(\vec{x}+\Delta\tau\vec{c}_{ij}) = g_{ij}(\vec{x}) + \Delta\tau \, \Omega_{ij}(\vec{x}) , \label{eq:LBM:residual2}
\end{equation}
where $\Omega_{ij}$ denotes the MRT collision operator for elasticity discussed in \cref{sec:LBM:elasticity}.
Hence, our goal is a scheme that rapidly yields $r_{ij} \to 0$.
This is achieved using the multi-grid method via 3 steps:
First, one applies a smoother to reduce the high-frequency errors, then one updates the solution with the so-called defect correction on a coarser grid, and finally applies the smoother again to remove high-frequency errors reintroduced by the update.
In the following, we first discuss the use of a under-relaxed LBM step as smoother, and afterwards the linear defect correction including the transfer between grids.
Note that we apply a standard linear multi-grid method as the considered elasticity equation is linear.

As realized in previous work, see e.g. \cite{Mavriplis2006,Patil2014}, the application of one LBM step is very similar to one step of the Jacobi iteration used for solving large equation systems.
In particular, the convergence of both methods deteriorates for elliptic problems of increasing size.
However, it is a suitable smoother when used as under-relaxation with parameter $\gamma \in (0,1)$, i.e.,
\begin{equation}
    g^{(k+1)}_{ij} = (1-\gamma) g^{(k)}_{ij} + \gamma g^{(k),\,**}_{ij} , \label{eq:MG:smoother}
\end{equation}
where we apply $K \in \N$ such under-relaxed LBM steps before and after the defect correction.
Hence, we start from $g^{(0)}_{ij} = g_{ij}(\tau)$, use \cref{eq:MG:smoother} for $k = 0,\dots,K-1$,
apply the defect correction to obtain $g^{(K+1)}_{ij}$, use \cref{eq:MG:smoother} for $k = K+1,\dots,2K$, 
and finally obtain $g_{ij}(\tau+\Delta\tau) = g^{(2K+1)}_{ij}$.

\begin{figure}[tbp]
    \centering
    \includegraphics{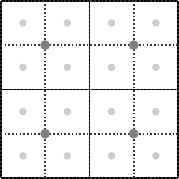}
    \caption{Sketch of a $4 \times 4$ fine grid (dotted lines) with cell centers (light gray circles) and the corresponding $2 \times 2$ coarse grid (solid lines) with cell centers (dark gray circles) as used in the multi-grid LBM.}
    \label{fig:grids}
\end{figure}

For the linear defect correction on a coarser grid, we simply double the lattice width $\Delta x$ to obtain the next coarser one.
As we use an cell-centered LBM, each cell of the coarse grid coincides with four cells of the finer one, see \cref{fig:grids}.
First, one computes the residual $r^{(K)}_{ij}$ on the current grid using \cref{eq:LBM:residual}.
To transfer the residual to the coarse grid, we use a simple (full) summation as restriction operator $\mathcal{R}$, i.e.,
the values of the four fine cells are added to obtain the value of the corresponding coarse cell, or in stencil notation:
\begin{equation}
    \mathcal{R} =
    \begin{bmatrix}
        1 & 1 \\ 1 & 1
    \end{bmatrix}.
    \label{eq:MG:restriction}
\end{equation}
This yields the coarse-grid residual $r^\text{c}_{ij} = \mathcal{R} r^{(K)}_{ij}$.
Vice versa, the coarse-correction values $g^\text{c}_{ij}$ are transferred to the current grid via the (full) prolongation $\mathcal{P} = \mathcal{R}^\top$, i.e., the correction value of all 4 fine cells is set to the value of the corresponding coarse cell.
To obtain the coarse-grid correction, one again has to solve \cref{eq:LBM:residual,eq:LBM:residual2} on the coarse grid, but with a different collision operator.
On the coarse grid, the contribution of the right-hand side force $\vec{f}_\text{eff}$ in the moment update \cref{eq:moments:halfForce} is removed up to the self-dependence $\vec{f}_*$ through $\nabla p$, viz.
\begin{equation}
    \vec{f}_*(g_{ij}) = \frac{\alpha^2}{2\varepsilon\Delta t}\nabla^c \big(\nabla\cdot\vec{\eta}(g_{ij})\big) ,\qquad
    \Tilde{\vec{f}}_\text{eff} = \vec{f}_\text{eff} - \vec{f}_*,
    \label{eq:MG:selfdependence}
\end{equation}
where $\Tilde{\vec{f}}_\text{eff}$ contains only explicit contributions, i.e., only depends on values of the previous time steps.
Moreover, the pre-collision distributions $g_{ij}$ must be reduced by the coarse-grid residual $r^\text{c}_{ij}$,
leading to \cref{alg:LBM:MG:residual,alg:LBM:MG:residual2,alg:LBM:MG:residual3} of \cref{alg:LBM:MG}.
Finally, the boundary conditions on the coarse grid are the homogeneous versions of those applied on the fine grid.
To solve for the coarse-grid correction, one can thus apply the multi-grid method recursively for $L \in \N$ levels with adapted arguments.

In summary, this leads to \cref{alg:LBM:MG} for the multi-grid LBM method for linear elasticity.
To include it into the overall poroelasticity LBM, lines 7 and 8 of \cref{alg:LBM} must be adapted, leading to \cref{alg:LBM:complete}.

\begin{algorithm}[tbp]
    \caption{LBM multi-grid step for quasi-static linear elasticity.}
    \label{alg:LBM:MG}
    \begin{algorithmic}[1]
        \Function{ELASTIC\_MULTIGRID}{$g_{ij}$, $r_{ij}$, $\Tilde{\vec{f}}_\text{eff}$, $\vec{f}_*$, $\ell$, $K$, $\gamma$}
            \LComment{Pre-smoothing}
            \For{$k = 1,\dots,K$}
                \State Compute $g^{**}_{ij} = \Call{ELASTIC\_STEP}{g_{ij} - r_{ij}, \Tilde{\vec{f}}_\text{eff} + \vec{f}_*(g_{ij})}$ \Comment{\cref{alg:LBM:elasticity}} \label{alg:LBM:MG:residual}
                \State Update $g_{ij} = (1-\gamma) g_{ij} + \gamma g^{**}_{ij}$                       \Comment{\cref{eq:MG:smoother}}
            \EndFor
            \LComment{Recursive defect correction}
            \If{$\ell > 1$}
                \State Restrict the residual $r_{ij}^\text{c} = \mathcal{R}\big(g_{ij} - \Call{ELASTIC\_STEP}{g_{ij} - r_{ij}, \Tilde{\vec{f}}_\text{eff} + \vec{f}_*(g_{ij})}\big)$ \Comment{\cref{eq:MG:restriction,alg:LBM:elasticity}} \label{alg:LBM:MG:residual2}
                \State Compute $g_{ij}^\text{c} = \Call{ELASTIC\_MULTIGRID}{\vec{0}, r_{ij}^\text{c}, \vec{0}, \vec{f}_*, \ell-1, K, \gamma}$
                \State Prolongate and update the solution $g_{ij} = g_{ij} + \mathcal{P}(g_{ij}^\text{c})$ \Comment{\cref{eq:MG:restriction}}
            \EndIf
            \LComment{Post-smoothing}
            \For{$k = 1,\dots,K$}
                \State Compute $g^{**}_{ij} = \Call{ELASTIC\_STEP}{g_{ij} - r_{ij}, \Tilde{\vec{f}}_\text{eff} + \vec{f}_*(g_{ij})}$ \Comment{\cref{alg:LBM:elasticity}} \label{alg:LBM:MG:residual3}
                \State Update $g_{ij} = (1-\gamma) g_{ij} + \gamma g^{**}_{ij}$                       \Comment{\cref{eq:MG:smoother}}
            \EndFor
            \State \Return $g_{ij}$
        \EndFunction
    \end{algorithmic}
\end{algorithm}

\begin{algorithm}[tbp]
    \caption{Semi-implicit LBM with multi-grid for poroelasticity.}
    \label{alg:LBM:complete}
    \begin{algorithmic}[1]
    \State Set $t = 0$ and initialize the distribution functions $f_i$ and $g_{ij}$                     \Comment{\cref{eq:ICs,eq:ICs:moments}}
    \For{$i = 1, \dots, N_T$}
        \LComment{Elasticity update}
        \State Set $\tau = 0$ and compute the explicit part of $s_\text{eff}^\text{im}(t)$             \Comment{\cref{eq:source:implicit}}
        \For{$j = 1, \dots, N_E$}
            \State Compute $s_\text{eff}^\text{im}(\tau,t)$ and $p(\tau,t)$         \Comment{\cref{eq:pressure,eq:source:implicit}}
            \State Compute the forces $\Tilde{\vec{f}}_\text{eff}(\tau,t)$ and $\vec{f}_*$ using $\nabla p(\tau,t)$ \Comment{\cref{eq:pressure:gradient,eq:MG:selfdependence}}
            \State Compute $g_{ij}(\tau+\Delta\tau,t) = \Call{ELASTIC\_MULTIGRID}{g_{ij}(\tau,t), \vec{0}, \Tilde{\vec{f}}_\text{eff}(\tau,t), \vec{f}_*, L, K, \gamma}$ \Comment{\cref{alg:LBM:MG}}
            \State Compute $\nabla\cdot\vec{\eta}(\tau+\Delta\tau,t)$                                       \Comment{\cref{eq:moments:bared,eq:LBM:divergence}}
            \State Advance pseudo time $\tau$ by $\Delta\tau$
        \EndFor
        \State Compute $\vec{\eta}(t)$ and $\tensor{\upsigma}(t)$                                                      \Comment{\cref{eq:displacement}}
        \LComment{Flow update}
        \State Compute $s_\text{eff}^\text{im}(t)$ and $p(\tau,t)$                         \Comment{\cref{eq:pressure,eq:source:implicit}}
        \State Compute $f_i(t+\Delta t) = \Call{FLOW\_STEP}{f_i(t), s_\text{eff}^\text{im}(t)}$ \Comment{\cref{alg:LBM:flow}}
        \State Advance time $t$ by $\Delta t$
    \EndFor
    \end{algorithmic}
\end{algorithm}

\subsection{Initial and boundary conditions}
\label{sec:LBM:BCs}

Next, we address the LBM implementation of the used initial and boundary conditions.
In particular, we discuss periodic, Dirichlet and Neumann boundary conditions for diffusive flow and elasticity.
Since all considered domains are axis-aligned and rectangular, we chose a cell-centered lattice, such that it is sufficient to use simple boundary schemes of (anti-)bounce-back type.
For more general boundary conditions and complex boundary shapes, we refer to \cite[Sec. 8.5]{Krueger2017} for diffusive flow and \cite{Boolakee2023BC} for linear elasticity.

\paragraph{Initial conditions}
The initial conditions for all numerical examples given below are
\begin{equation*}
    p \big|_{t = 0} = 0, \qquad \vec{\eta} \big|_{t = 0} = \vec{0} .
\end{equation*}
However, to maintain second-order convergence, the initial conditions in the LBM must be chosen such that the physical initial conditions are satisfied in \cref{eq:pressure,eq:displacement} at $t = 0$ \cite{Dellar2013,Krueger2017}.
Hence, the initial distribution functions for pressure and displacement are
\begin{equation}
    f_i = - \frac{\Delta t}{2} w_i \big( s - \kappa \nabla\cdot\vec{g} \big) \big|_{t = 0} ,\qquad
    g_{ij} = \sum_\nu c_{ij}^\nu m^0_\nu .\label{eq:ICs}
\end{equation}
Here, the effective source term at time $t = 0$ disregards the unknown volumetric source $-\frac{\alpha}{\varepsilon} \partial_{t}\nabla\cdot\vec{\eta}$.
Moreover, the initial moments  at time $t = 0$ are the force-corrected equilibrium moments
\begin{equation}
    \begin{pmatrix} m^0_{10} \\[.5ex] m^0_{01} \end{pmatrix} = -\frac{\varepsilon^2}{2} \vec{f}_\text{eff} \big|_{t = 0}, \qquad
    m^0_{s} = m^0_{d} = m^0_{11} = 0, \qquad
    \begin{pmatrix} m^0_{12} \\[.5ex] m^0_{21} \end{pmatrix} = -\frac{\varepsilon^2}{6} \vec{f}_\text{eff} \big|_{t = 0} ,\qquad
    m^0_{22} = 0 ,\label{eq:ICs:moments}
\end{equation}
where the effective force at time $t = 0$ contains only known terms since $\vec{f}_\text{eff} |_{t = 0} = ( \vec{f} - \alpha \nabla p ) |_{t = 0}$.
Note that instead of explicitly initializing the distribution functions $g_{ij}$ for elasticity, one could just use more iterations in the first time step to achieve the stationary solution.

\paragraph{Periodic boundary conditions}

Periodic boundary conditions on an open and connected subset of the boundary $\Gamma_\text{\!per} \subset \partial\mathcal{D}$ describe solutions repeating infinitely in space by using the equality of the solution at opposite boundaries of the finite domain $\mathcal{D} = (0, a) \times (0, b)$, $a,b \in \mathbb{R}$, i.e.,
\begin{equation*}
    p \big|_{x_1 = 0} = p \big|_{x_1 = a} , \qquad
    p \big|_{x_2 = 0} = p \big|_{x_2 = b} , \qquad
    \vec{\eta} \big|_{x_1 = 0} = \vec{\eta} \big|_{x_1 = a} , \qquad
    \vec{\eta} \big|_{x_2 = 0} = \vec{\eta} \big|_{x_2 = b} , \qquad
    \text{on } (0, t_f) \times \Gamma_\text{\!per}.
\end{equation*}
This requires that distribution functions leaving the domain on one side, simultaneously, re-enter at the opposite side.
Since the lattice is cell-centered, this can be realized by simply copying the outgoing distribution functions during streaming to the respective cells on the opposite side.
Furthermore, this scheme maintains the convergence order of the underlying LBM.
Note that this scheme is automatically included in the streaming step, when a rolling shift implementation is used, e.g., by the function \texttt{roll} in the libraries \texttt{numpy} or \texttt{jax} of Python.

\begin{figure}[tb]
    \centering
    \includegraphics{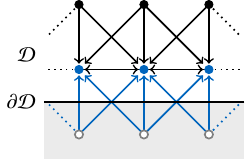}
    \caption{Missing incoming distribution functions (blue) at the bottom boundary of the rectangular domain $\mathcal{D}$ and known incoming distribution functions (black) from inside the domain.
        Black dots: interior nodes, blue dots: boundary nodes, empty gray nodes: ghost nodes.}\label{fig:boundary}
\end{figure}

In case of other boundary conditions than periodic ones, there are missing distribution functions after the streaming which would come from outside the domain to the boundary nodes, see \cref{fig:boundary} for a sketch at the bottom boundary.
These missing distributions can be reconstructed from the post-collision distribution functions taking into account the respective boundary condition.
To this end, the basic idea of bounce-back schemes is that distribution functions hitting the boundary are reflected, and thus the missing distributions after streaming are computed from the leaving post-collision distributions in opposite directions.

\paragraph{Dirichlet boundary conditions}

First, we consider a Dirichlet boundary condition for the pressure $p$ to be prescribed on an open and connected subset of the boundary $\Gamma_{D,p} \subset \partial\mathcal{D}$, viz,
\begin{equation*}
    p = p_D \qquad\text{on } (0, t_f) \times \Gamma_{D,p} ,
\end{equation*}
for some given continuous function $p_D : (0, t_f) \times \Gamma_{D,p} \to \mathbb{R}$.
In this case, we apply the second-order consistent, halfway anti-bounce-back scheme following \cite{Ginzburg2005}.
When a boundary node with position $\vec{x}$ connects via the velocity direction $\vec{c}_i$ to the Dirichlet boundary $\Gamma_{D,p}$, then the distribution function at the next time step is given by
\begin{equation}
    f_{\check i}(t+\Delta t, \vec{x}) = -f^*_i(t, \vec{x}) + 2 w_i p_D\big(t+\Delta t, \vec{x}+\tfrac12\Delta t\vec{c}_i\big) .\label{eq:BC:pressure}
\end{equation}
Here, $\check i$ denotes the index of the opposite direction of index $i$, e.g. $\check 4 = 2$ and $\check 5 = 7$, cf. \cref{fig:velocities}.

Similarly, when a Dirichlet boundary condition for the displacement $\vec{\eta}$ at an open and connected subset of the boundary $\Gamma_{D,\vec{\eta}} \subset \partial\mathcal{D}$ is given by
\begin{equation*}
    \vec{\eta} = \vec{\eta}_D \qquad\text{on } (0, t_f) \times \Gamma_{D,\vec{\eta}} ,
\end{equation*}
for some given continuous function $\vec{\eta}_D : (0, t_f) \times \Gamma_{D,\vec{\eta}} \to \mathbb{R}$, then we use the fullway bounce-back scheme developed in \cite{Boolakee2023BC}.
For a boundary node with position $\vec{x}$ connecting to the Dirichlet boundary $\Gamma_{D,\vec{\eta}}$ via the velocity direction $\vec{c}_{ij}$, the distribution function at the next pseudo-time step is given by
\begin{equation}
    g_{\bar i \bar j}(\tau+\Delta\tau, t, \vec{x}) = g_{ij}(\tau-\Delta\tau, t, \vec{x}) - 6 w_{ij} \vec{c}_{ij} \cdot \vec{\eta}_D\big(t, \vec{x}+\tfrac12\Delta t\vec{c}_{ij}\big) .\label{eq:BC:elasticity}
\end{equation}
As before, we use the 2D Miller index notation here, i.e., $\bar 1 = -1$, cf. \cref{fig:velocities}.
Note that the delay of one pseudo-time step is equivalent to two full streaming steps into the respective ghost node and back, hence the term ``fullway''.
Furthermore, this delay is necessary to suppress grid-scale oscillations in the displacement solution according to \cite{Boolakee2023BC}.
As shown therein, this boundary scheme yields a second-order consistent displacement solution and a first-order consistent stress solution for pseudo time $\tau \to \infty$.

\paragraph{Neumann boundary conditions}

For the pressure $p$, we consider homogeneous Neumann boundary conditions (no-flow conditions) posed on an open and connected subset of the boundary $\Gamma_{N,p} \subset \partial\mathcal{D}$, viz.
\begin{equation*}
    \kappa \nabla p \cdot \vec{n} = 0 \qquad\text{on } (0, t_f) \times \Gamma_{N,p} ,
\end{equation*}
where $\vec{n}$ denotes the outward unit normal at the boundary $\Gamma_{N,p}$.
This condition can be implemented using the simple halfway bounce-back scheme \cite[Sec. 8.5.3]{Krueger2017}.
The new distribution function at a boundary node connecting via the velocity direction $\vec{c}_i$ to the Neumann boundary $\Gamma_{N,p}$ is thus given by
\begin{equation}
    f_{\check i}(t+\Delta t, \vec{x}) = f^*_i(t, \vec{x}) ,\label{eq:BC:noFlow}
\end{equation}
where $\check i$ denotes the opposite direction index introduced above for Dirichlet conditions.
Note that this boundary scheme is first-order consistent in general, but second-order consistent if the flux tangential to the boundary vanishes as considered in the numerical examples below.

Finally, a Neumann boundary condition for the displacement $\vec{\eta}$ (traction condition) can be imposed at an open and connected subset of the boundary $\Gamma_{N,\vec{\eta}} \subset \partial\mathcal{D}$, i.e.,
\begin{equation*}
    \tensor{\upsigma} \vec{n} = \vec{t} \qquad\text{on } (0, t_f) \times \Gamma_{N,\vec{\eta}} ,
\end{equation*}
with the outward unit normal $\vec{n}$ and some given continuous traction function $\vec{t} : (0, t_f) \times \Gamma_{N,\vec{\eta}} \to \mathbb{R}^2$.
Following \cite{Boolakee2023BC}, we apply a second-order consistent, halfway anti-bounce-back scheme to obtain the missing distribution function at a boundary node connecting via the velocity direction $\vec{c}_i$ to the Neumann boundary $\Gamma_{N,\vec{\eta}}$.
This yields
\begin{equation}
    g_{\bar i\bar j}(\tau+\Delta\tau, t, \vec{x}) = -g_{ij}(\tau, t, \vec{x}) + \varepsilon\, \varsigma_{\bar i\bar j}\big(t, \vec{x}+\tfrac12\Delta t\vec{c}_{ij}\big) ,\label{eq:BC:traction}
\end{equation}
where the traction correction terms $\varsigma_{ij}$ are obtained via the stress \cref{eq:displacement} from the traction $\vec{t}$ and the normal direction $\vec{n}$.
For the left boundary, where $\vec{n} = (-1, 0)^\top$, one gets
\begin{equation*}
    \varsigma_{10} = -\upsigma_{11} = t_1, \qquad \varsigma_{11} = -\frac12 \upsigma_{12} = \frac12 t_2, \qquad \varsigma_{1\bar 1} = \frac12 \upsigma_{12} = -\frac12 t_2 .
\end{equation*}
Analogous formulas hold for the right, bottom and top boundaries, but are not explicitly shown here.

\section{Numerical Experiments}
\label{sec:numerics}

This section presents several numerical examples to study the numerical accuracy and stability.
To this end, we consider three problem settings: the first one is a periodic problem with manufactured solution to assess the performance of the proposed LBM in the bulk without interference of boundary, while the second problem is the two-dimensional version of the classical (quasi-1D) Terzaghi consolidation problem, which even has a discontinuous solution at the initial time due to the instantaneous response to a traction boundary.
Finally, the third example extends the latter setting to a truly two-dimensional problem.
We present detailed parameter studies concerning the coupling strength, i.e., the Biot--Willis coefficient $\alpha \in (0,1]$, as well as the convergence with respect to the lattice width $\Delta x$ and the number of pseudo-time steps $N_E$, while the time step $\Delta t = \Delta x^2$ and the length ratio $\varepsilon = \Delta x$ follow the typical viscous and elastic scaling, respectively.
When applying $N_E$ steps of the multi-grid LBM of \cref{sec:LBM:MG} using $K$ smoothing steps with under-relaxation rate $\gamma = 3/4$, we fix the maximal grid level $L = \log_2(N_x)$, and abbreviate this as MG($N_E,K$) for simplicity.
Our implementation is based on XLB \cite{Ataei2024}, a \texttt{jax}-based Python library allowing for parallel LBM simulations.
The code is available in the GitHub repository \url{https://github.com/s-lunowa/xlb-biot}.
All presented simulations were run on an ordinary workstation with Intel i7-11700 CPU (16 $\times$ 2.50GHz).

For the convergence studies, the errors of the numerical solutions of the poroelasticity \cref{eq:elasticity:dim,eq:stress:dim,eq:flow:dim,eq:darcy:dim} at the times $t_k = k \Delta t$, $k \in \Set{0, 1, \dots, N_T}$, are measured based on the discrete $L^2(\mathcal{D})$-norm, i.e.,
\begin{align*}
    e^k_{p} &= \Norm*{ p(t_k) - p^\text{LBM}(t_k) }_{L^2(\mathcal{D})} = \Delta x \left( \sum_{m=1}^{N_x N_y} \left|p(t_k, \vec{x}_m) - p^\text{LBM}(t_k, \vec{x}_m)\right|^2 \right)^{1/2}, \\
    e^k_{\vec{\eta}} &= \Norm*{ \vec{\eta}(t_k) - \vec{\eta}^\text{LBM}(t_k) }_{L^2(\mathcal{D})} = \left( \sum_{i=1}^2 \Norm*{ \eta_i(t_k) - \eta^\text{LBM}_i(t_k) }_{L^2(\mathcal{D})}^2 \right)^{1/2}, \\
    e^k_{\tensor{\upsigma}} &= \Norm*{ \tensor{\upsigma}(t_k) - \tensor{\upsigma}^\text{LBM}(t_k) }_{L^2(\mathcal{D})} = \left( \sum_{ab \in \Set {11, 12, 22}} \Norm*{ \upsigma_{ab}(t_k) - \upsigma^\text{LBM}_{ab}(t_k) }_{L^2(\mathcal{D})}^2 \right)^{1/2} ,
\end{align*}
where $N_x$ and $N_y$ denote the number of lattice cells in $x_1$- and $x_2$-direction, respectively.
Furthermore, the time-space $L^2(\mathcal{Q})$-error of the pressure is given by
\begin{equation*}
    e_{p} = \Norm*{ p - p^\text{LBM} }_{L^2(\mathcal{Q})} = \left( \Delta t \sum_{k=0}^{N_T} \Norm*{ e^k_p }_{L^2(\mathcal{D})}^2 \right)^{1/2},
\end{equation*}
and analogously for the time-space $L^2(\mathcal{Q})$-error of displacement $e_{\vec{\eta}}$ and of stress $e_{\tensor{\upsigma}}$.
To have scaling-independent, easily interpretable results, we will primarily report relative errors such as $e^\text{rel}_{p} = e_{p} / \Norm{p}_{L^2(\mathcal{Q})}$.

\subsection{Periodic problem with manufactured analytical solution}
\label{sec:numerics:manufactured}

First, a manufactured solution to \cref{eq:elasticity:dim,eq:stress:dim,eq:flow:dim,eq:darcy:dim} is considered for numerical validation of accuracy and stability.
To avoid the influence of boundary conditions, a periodic problem ($\Gamma_\text{\!per} = \partial\mathcal{D}$) is posed in the domain $\mathcal{D} = (0, 1)^2$ up to the final time $t_f = 1$.
Inspired by \cite{Boolakee2023}, the manufactured solution used here is given by
\begin{align*}
    \vec{\eta}(t, \vec{x}) &= \frac{1 - \exp(-8\pi^2 \kappa t)}{2} \begin{pmatrix} 9 \cos(2\pi x_1) \sin(2\pi x_2) \\ 7 \sin(2\pi x_1) \cos(2\pi x_2) \end{pmatrix}, \\
    p(t, \vec{x}) &= - \frac{(16\lambda + 30\mu) \pi}{\alpha} \left(1 - \exp(-8\pi^2 \kappa t) \right) \sin(2\pi x_1) \sin(2\pi x_2) ,
\end{align*}
which leads to the following forces and source term
\begin{gather*}
    \vec{f}(t, \vec{x}) = \frac{1 - \exp(-8\pi^2 \kappa t)}{2} \begin{pmatrix} 16 \pi^2 \mu \cos(2\pi x_1) \sin(2\pi x_2) \\ 0 \end{pmatrix}, \qquad
    \vec{g}(t, \vec{x}) = \vec{0} ,\\
    s(t, \vec{x}) = -8 \pi^3 \kappa \left(16 \alpha \exp(-8\pi^2 \kappa t) + \frac{16 \lambda + 30 \mu}{\alpha}\right) \sin(2\pi x_1) \sin(2\pi x_2) .
\end{gather*}
Setting the (dimensionless) Young's modulus $E = 0.11$ and Poisson's ratio $\nu = 0.8$ as in \cite{Boolakee2023} (with plane-stress assumption), the physical parameters are
\begin{equation*}
    \lambda = \frac{E \nu}{1 - \nu^2} = \frac{11}{45} ,\qquad
    \mu = \frac{E}{2(1 + \nu)} = \frac{11}{360} ,\qquad
    \kappa = \frac{1}{10},
\end{equation*}
while we consider varying Biot--Willis coefficients $\alpha \in [0.5, 1.0]$.
In the following, we discuss the influence of this coefficient $\alpha$, of the number of pseudo-time steps $N_E$ and of the ratio $r \in \Set{0, 1/2, 1}$ in the semi-implicit discretization.

\paragraph{Weak coupling}

\begin{figure}[tbp]%
    \centering
    \includegraphics{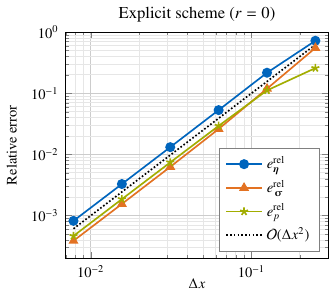}
    \includegraphics{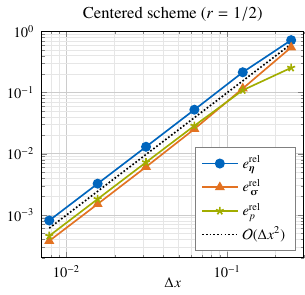}
    \includegraphics{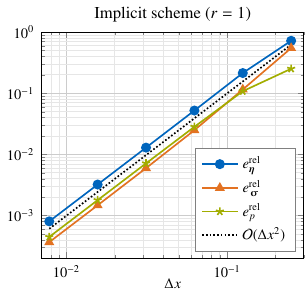}
    \caption{Convergence results for the periodic problem with $\alpha = 0.5$ using the explicit, centered and implicit coupling schemes with MG(3,2).}\label{fig:periodic:alpha05ne3}
\end{figure}

\begin{figure}[p]%
    \centering
    \includegraphics{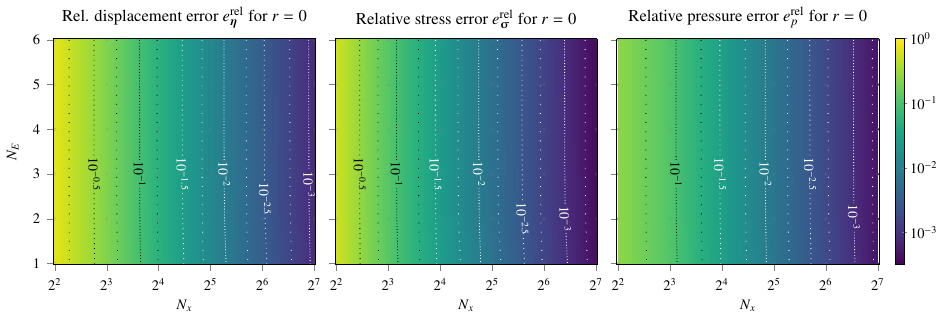}
    \smallskip

    \includegraphics{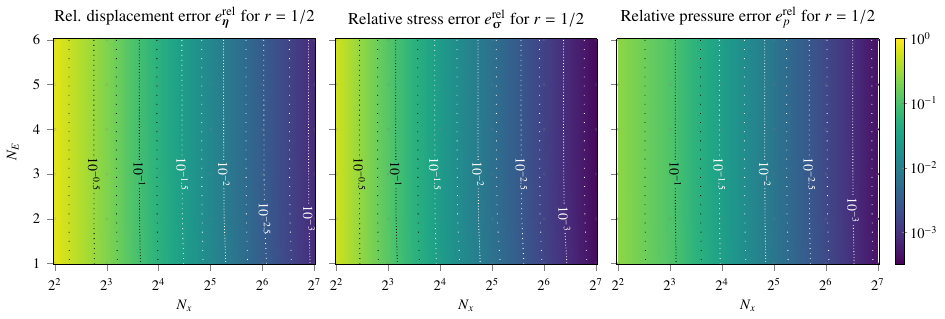}
    \smallskip

    \includegraphics{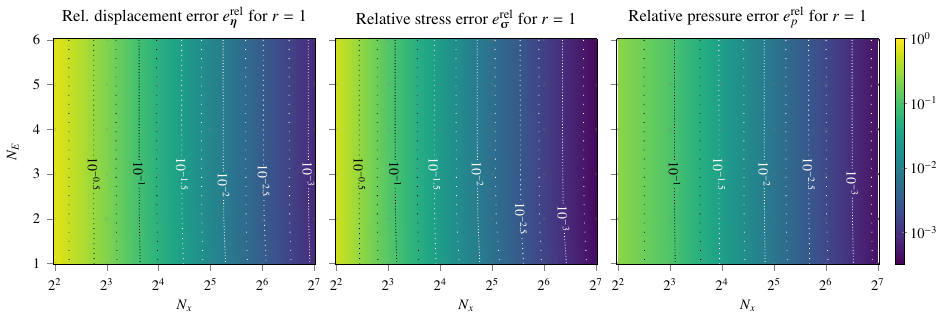}
    \caption{Convergence results for the periodic problem with $\alpha = 0.5$ using the explicit ($r=0$), centered ($r=1/2$) and implicit ($r=1$) coupling schemes (top to bottom) with MG($N_E$,2).}\label{fig:periodic:alpha05}
\end{figure}

First, we consider rather low Biot--Willis coefficients $\alpha \leq 0.5$, which means that the poroelastic system is weakly coupled.
In the following, we exemplary present the case $\alpha = 0.5$, but observe the same behavior also in the cases $\alpha \in \Set{0.3, 0.4}$.
Using a fixed number $N_E$ of multi-grid steps MG($N_E$,2), the algorithmic the algorithmic complexity of $\mathcal{O}(N_E N_x^2)$ per time step is optimal.
For $N_E = 3$, we observe second-order convergence in $\varepsilon = \Delta x = 1/N_x$ for the explicit ($r=0$), centered ($r=1/2$) and implicit ($r=1$) coupling scheme, see \cref{fig:periodic:alpha05ne3}.
In particular, we observe that the errors in pressure and stress are similar, while the displacement errors are about two times larger than the others.
A detailed convergence study for varying number $N_E$ of multi-grid steps is displayed in \cref{fig:periodic:alpha05} for $\alpha = 0.5$ and all three schemes.
While the numerical errors might be dominated by the elasticity approximation for $N_E \leq 2$ multi-grid steps, as visible by the slightly slanted contour lines, they follow the typical second-order convergence rate for $N_E \geq 3$, where the contour lines are vertically aligned.
In general, there are virtually no differences between the results of all three coupling schemes.

\begin{figure}[tbp]%
    \centering
    \includegraphics{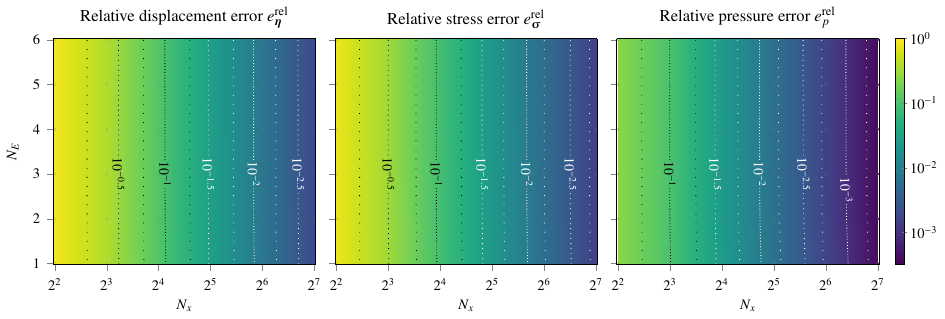}
    \caption{Convergence results for the periodic problem with $\alpha = 0.5$ using the simpler pressure gradient formulation of \cref{eq:pressure:gradient:alternative} in the centered coupling scheme ($r = 1/2$) with MG($N_E$,2).}\label{fig:periodic:alpha05_altgrad}
\end{figure}

Using the simpler pressure gradient formulation of \cref{eq:pressure:gradient:alternative}, we also observe second-order convergence, but an about two times increased elasticity error (in displacement and stress), as depicted in \cref{fig:periodic:alpha05_altgrad} for the centered scheme (all three schemes produce practically indistinguishable results).
Note however, that in case of strong coupling discussed below, the onset of instability for the explicit coupling scheme ($r=0$) occurs slightly later with that formulation, while the implicit coupling scheme ($r=1$) remains stable.

\begin{figure}[tbp]%
    \centering
    \includegraphics{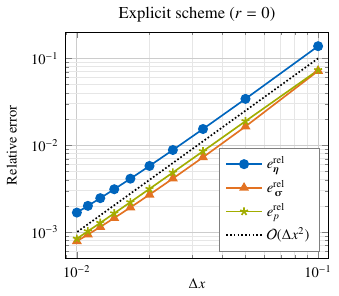}
    \hspace{-1em}
    \includegraphics{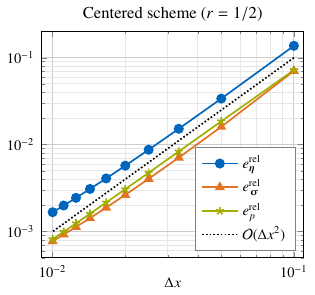}
    \hspace{-1em}
    \includegraphics{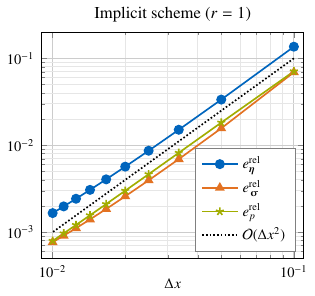}
    \caption{Convergence results for the periodic problem with $\alpha = 0.5$ using the explicit, centered and implicit coupling schemes with $N_E = 1.5N_x$ pseudo-time steps (no multi-grid).}\label{fig:periodic:alpha05ne15nx}
\end{figure}

\begin{figure}[tbp]%
    \centering
    \includegraphics{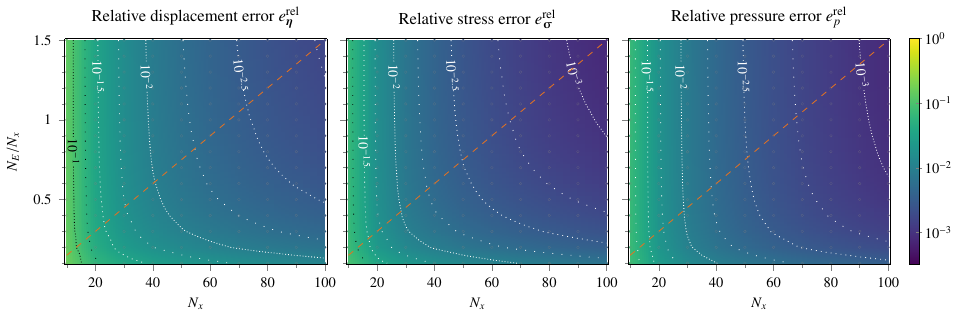}
    \caption{Convergence results for the periodic problem with $\alpha = 0.5$ using the centered coupling scheme ($r = 1/2$) with $N_E$ pseudo-time steps (no multi-grid).
        The orange dashed line $N_E = 0.015 N_x^2$ indicates the transition between dominant error due to spatial discretization (above the line) and due to insufficient pseudo-timestepping (below the line).}\label{fig:periodic:alpha05_stepping}
\end{figure}

Furthermore, we compare this to the simple \cref{alg:LBM} using pseudo-timestepping without multi-grid method for the elasticity.
Fixing the number of pseudo-time steps to $N_E = 1.5N_x$, the numerical results indicate almost second-order convergence, as depicted exemplary for $\alpha = 0.5$ in \cref{fig:periodic:alpha05ne15nx} for all three coupling schemes.
This linearly growing number of pseudo-time steps is not sufficient, as a detailed convergence study for varying $N_E$ shows, cf.~\cref{fig:periodic:alpha05_stepping} for the centered scheme.
While the numerical errors are dominated by the elasticity approximation for small $N_E$ as visible by the rather horizontal contour lines, they follow the typical second-order convergence rate for large $N_E$, where the contour lines are vertically aligned.
In particular, the transition between these two regimes seems to occur around $N_E \approx 0.015 N_x^2$, explaining the degradation of the convergence order visible for small $\Delta x$ in \cref{fig:periodic:alpha05ne15nx}.
Note that this quadratic scaling is characteristic for point-wise iterative solvers when applied to second-order elliptic problems.
Hence, the multi-grid method achieves second-order convergence at significantly reduced computational cost.

\paragraph{Strong coupling}

\begin{figure}[tbp]%
    \centering
    \includegraphics{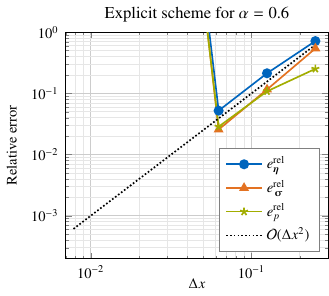}
    \includegraphics{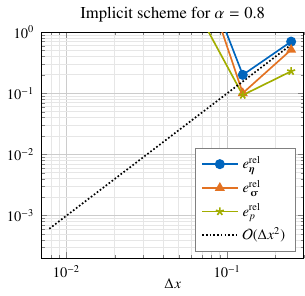}
    \includegraphics{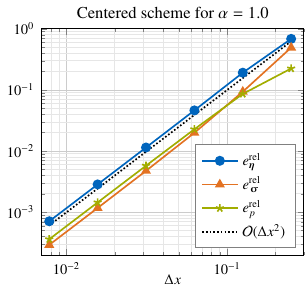}
    \caption{Divergence results for the periodic problem using the explicit coupling scheme ($r = 0$) for $\alpha = 0.6$ and the implicit coupling scheme ($r = 1$) for $\alpha = 0.8$, but convergence for the centered coupling scheme ($r = 1/2$) for $\alpha = 1.0$, each using MG(4,2).}\label{fig:periodic:divergence}
\end{figure}

Next, we consider large Biot--Willis coefficients $\alpha \in \Set{0.6, 0.7, 0.8, 0.9, 1.0}$, which means that the poroelastic system is strongly coupled.
Here, both the explicit scheme and (slightly later) the implicit scheme show instabilities for a sufficiently large number of multi-grid steps $N_E$, which decreases with increasing Biot--Willis coefficient $\alpha$.
In particular, when using $N_E = 4$ multi-grid steps, we observe divergence of the errors for the explicit scheme when $\alpha \geq 0.6$, and for the implicit scheme when $\alpha \geq 0.8$, as depicted in \cref{fig:periodic:divergence}.
Hence, the implicit scheme is slightly better than the explicit one, but ultimately also unstable for $\alpha$ close to 1.
Similar behavior has also been observed for other numerical methods due to spurious oscillations between stress and pressure, cf. \cite{Phillips2008,Bean2014,Terekhov2020}.

\begin{figure}[tbp]
    \includegraphics{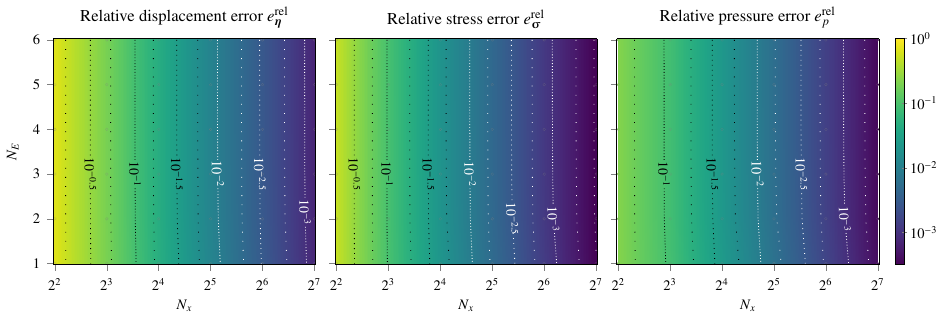}
    \caption{Convergence results for the periodic problem with $\alpha = 1.0$ using the centered coupling scheme ($r = 1/2$) with MG($N_E$,2).}\label{fig:periodic:alpha10}
\end{figure}

On the other hand, the centered scheme still gives second-order convergence for all $\alpha$ up to one, see \cref{fig:periodic:divergence} for $\alpha = 1$.
Probably, this is due to averaging over two time step which leads to the cancellation of those spurious oscillations.
Furthermore, a numerical convergence study varying the number $N_E$ of multi-grid steps yields again a clear second-order convergence for the centered scheme when $N_E \geq 3$, cf.~\cref{fig:periodic:alpha10} for $\alpha = 1.0$.
Note that we observe analogous results for $\alpha \in \Set{0.6, 0.7, 0.8, 0.9}$.

\subsection{Terzaghi's consolidation problem}
\label{sec:numerics:terzaghi}

The second poroelastic example originates from Terzaghi's consolidation problem \cite{Terzaghi1923}, which describes a fluid-filled soil layer which becomes loaded from above.
This leads to a so-called overpressurization of the fluid which continually vanishes due to flow through the surface.
Consequently, the elastic soil bears more and more of the loading and subsides.
This problem has an analytical solution and, hence, has been considered as benchmark for many numerical methods, see e.g. \cite{Murad1994,Phillips2008,Choo2018,Kadeethum2021} for finite elements, \cite{Honorio2018,Sokolova2019} for finite volumes,
and \cite{Zhao2013} for a direct numerical simulations of the averaged Navier–Stokes equation coupled to discrete particles.

While the original problem considers the quasi-1D solution to the three-dimensional poroelastic equations ($d = 3$), we adapt the quasi-1D formulation for $d = 2$.
Following the derivation in \cite[Secs.~5.1 \& 5.2.2]{Coussy2004}, we consider \cref{eq:elasticity:dim,eq:stress:dim,eq:flow:dim,eq:darcy:dim} in a domain $\mathcal{D} = (0, l) \times (0,1)$ until the final time $t_f = 1$ with fixed pressure and traction at the top
\begin{equation*}
    p = 0 \qquad\text{and}\qquad
    \tensor{\upsigma} \vec{n} = \begin{pmatrix} 0 \\ -\varpi \end{pmatrix} \qquad\text{on } (0,t_f) \times \Set{ \vec{x} \in \partial\mathcal{D} \given x_2 = 1 },
\end{equation*}
while no-flow condition and fixed displacement are prescribed at the bottom
\begin{equation*}
    \kappa \nabla p \cdot \vec{n} = 0 \qquad\text{and}\qquad
    \vec{\eta} = \vec{0} \qquad\text{on } (0,t_f) \times \Set{ \vec{x} \in \partial\mathcal{D} \given x_2 = 0 },
\end{equation*}
and periodic boundary conditions are posed on the left and right at $x_1 = 0$ and $x_1 = l$.
Note that the length $l > 0$ is irrelevant since the solution is independent of the first spatial coordinate and we choose $l = 1$ in the numerical experiments.
No external forces or source are present, i.e., $\vec{f} \equiv \vec{g} \equiv \vec{0}$ and $s \equiv 0$.
Moreover, the initial conditions are given by
\begin{equation*}
    p \big|_{t = 0} = 0 \qquad\text{and}\qquad
    \vec{\eta} \big|_{t = 0} = \vec{0} \qquad\text{in } \mathcal{D}.
\end{equation*}
Due to the incompatibility of the initial and boundary conditions, the solution has a discontinuity at time $t = 0$ with an instantaneous response
\begin{equation*}
    p(0^+, \vec{x}) = p_0 = \frac{\alpha \varpi}{\lambda+2\mu + \alpha^2} \qquad\text{and}\qquad
    \vec{\eta}(0^+, \vec{x}) = \begin{pmatrix} 0 \\ - p_0 x_2 / \alpha \end{pmatrix} \qquad\text{for }\vec{x} \in \mathcal{D}.
\end{equation*}
Hereafter, the solution evolves continuously.
In particular, the pressure and the surface subsidence $S = -\eta_2 |_{x_2 = 1}$ are explicitly given by
\begin{align*}
    p(t, \vec{x}) &= p_0 \sum_{k = 0}^\infty \frac{2 (-1)^k}{\big(k + \tfrac12\big)\pi} \exp\Big(-\big(k + \tfrac12\big)^2 \pi^2 c_f t\Big) \cos\left(\big(k + \tfrac12\big) \pi x_2 \right) &&\text{for }(t, \vec{x}) \in (0, \infty) \times \mathcal{D}, \\
    S(t) &= S_\infty + (S_0 - S_\infty) \sum_{k = 0}^\infty \frac{2}{\big(k + \tfrac12\big)^2 \pi^2} \exp\Big(-\big(k + \tfrac12\big)^2 \pi^2 c_f t\Big) &&\text{for }t \in (0, \infty),
\end{align*}
where the effective diffusivity $c_f$ as well as the initial and final surface subsidence are defined as
\begin{equation*}
    c_f = \frac{\lambda + 2 \mu}{\lambda + 2 \mu + \alpha^2} \kappa ,\qquad
    S_0 = \frac{\varpi}{\lambda + 2 \mu + \alpha^2} ,\qquad
    S_\infty = \frac{\varpi}{\lambda + 2 \mu}.
\end{equation*}
Setting the (dimensionless) Young's modulus $E = 1$ and Poisson's ratio $\nu = 0.8$ (with plane-stress assumption), we choose the physical parameters
\begin{equation*}
    \lambda = \frac{E \nu}{1 - \nu^2} = \frac{20}{9} ,\qquad
    \mu = \frac{E}{2(1 + \nu)} = \frac{5}{18} ,\qquad
    \kappa = \frac{\lambda + 2 \mu + \alpha^2}{ \lambda + 2 \mu } = \frac{25 + 9\alpha^2}{25}, \qquad
    \varpi = \alpha + \frac{\lambda + 2 \mu}{\alpha} ,
\end{equation*}
which implies $p_0 = 1$.
Note that due to the permeability $\kappa > 1$, we use a scaling $\Delta t = \Delta x^2 / 4$ to reduce the numerical errors in the pressure solution.
Since only the pressure and the surface subsidence are given analytically, we report relative errors of pressure $e^\text{rel}_{p}$ and of subsidence $e^\text{rel}_{S} = \Norm{S - S^\text{LBM}}_{L^2(0,T)} / \Norm{S}_{L^2(0,T)}$, where the numerical subsidence is defined as average of $-\eta_2$ over the last row of cells next to the boundary at $x_2 = 1$.

\begin{figure}[tbp]
    \centering
    \includegraphics{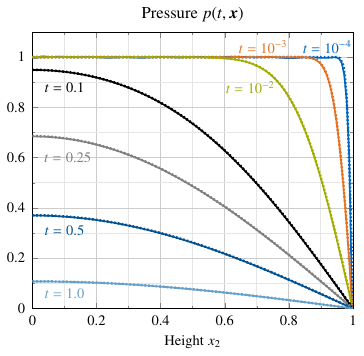}
    \hspace{2em}
    \includegraphics{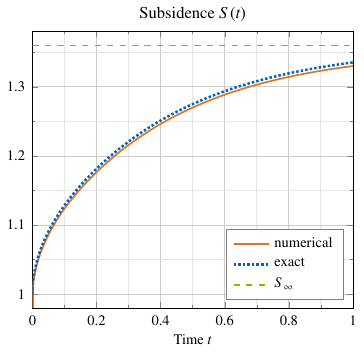}
    \caption{Exact and numerical solutions for Terzaghi's problem with $\alpha = 1.0$ using the centered scheme ($r = 1/2$) with MG(3,4) for $N_x = 128$. In the left sub-figure, dotted lines represent the exact solution, while solid lines depict the numerical solution.}\label{fig:Terzaghi:solution}
\end{figure}

\begin{figure}[tbp]
    \centering
    \includegraphics{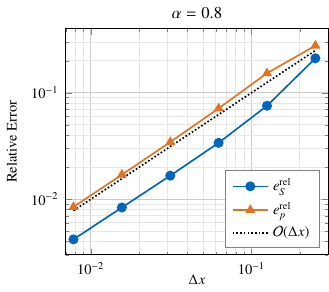}
    \hspace{-1ex}
    \includegraphics{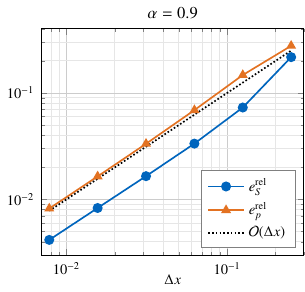}
    \hspace{-1ex}
    \includegraphics{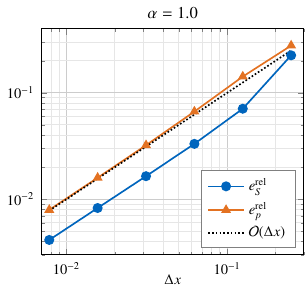}
    \caption{Convergence results for Terzaghi's problem with $\alpha = 0.8$, $0.9$ and $1.0$ using the centered coupling scheme ($r = 1/2$) with MG(3,4).}\label{fig:Terzaghi:centered}
\end{figure}

\begin{figure}[tbp]
    \centering
    \includegraphics{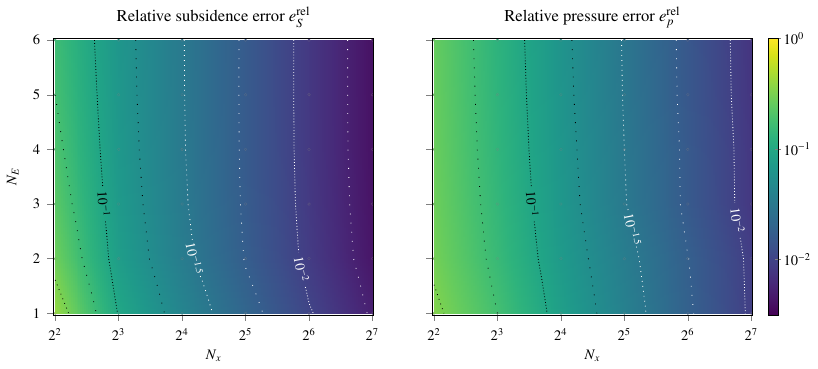}
    \caption{Convergence results for Terzaghi's problem with $\alpha = 1.0$ using the centered coupling scheme ($r = 1/2$) with MG($N_E$,4).}\label{fig:Terzaghi:alpha10}
\end{figure}

The numerical solution of the centered coupling scheme ($r = 1/2$) for $N_x = 128$ using MG(3,4) and the analytical solution are depicted in detail for $\alpha = 1$ in \cref{fig:Terzaghi:solution}.
The results excellently match already at time $t = 10^{-4}$ (after only 8 time steps), thus illustrating that the LBM scheme is well able to reproduce the solution due to instantaneous loading.
Particularly, we observe first-order convergence of the centered coupling scheme with MG(3,4) for varying Biot--Willis coefficient $\alpha \in \Set{ 0.8, 0.9, 1.0}$, as shown in \cref{fig:Terzaghi:centered}.
Note that convergence order is reduced due to the implementation of the boundary conditions and the discontinuity at time $t = 0$.
In contrast to the previous numerical example, also the explicit scheme ($r = 0$) and the implicit scheme ($r = 1$) feature such a convergence behavior for all $\alpha \in \Set{0.8, 0.9, 1.0}$.
Concerning the influence of the number of multi-grid steps $N_E$, a convergence study for the centered coupling scheme ($r = 1/2$) is displayed in \cref{fig:Terzaghi:alpha10} for $\alpha = 1$.
Analogously to the periodic example, the numerical errors are dominated by the elasticity approximation for small
$N_E \leq 2$, while they are dominated by the first-order spatio-temporal error for larger $N_E$.
Here however, the transition region seems slightly wider.
Note that we observe very similar results for $\alpha \in \Set{0.8, 0.9}$ and all three coupling schemes.

\subsection{Two-dimensional loading problem}
\label{sec:numerics:2D}

The third example extends Terzaghi's consolidation problem to two dimensions by introducing an position-dependent load.
Again, this describes a fluid-filled soil layer which becomes locally loaded from above.
The resulting overpressurized fluid slowly flows out through the surface, so that the elastic soil bears more and more of the loading and subsides.
While Terzaghi's problem considers a quasi-1D solution, here we also observe horizontal displacement due to the nonuniform load.

To this end, we consider the dimensional \cref{eq:elasticity,eq:stress,eq:flow,eq:darcy} in the domain $\mathcal{D} = (0, 100)^2$ until the final time $t_f = 2 \cdot 10^6$.
No external forces or source are present, i.e., $\vec{f} \equiv \vec{g} \equiv \vec{0}$ and $s \equiv 0$, and the initial conditions are
\begin{equation*}
    p \big|_{t = 0} = 0 \qquad\text{and}\qquad
    \vec{\eta} \big|_{t = 0} = \vec{0} \qquad\text{in } \mathcal{D}.
\end{equation*}
At the top boundary, pressure and traction are imposed
\begin{equation*}
    p = 0 \qquad\text{and}\qquad
    \tensor{\upsigma} \vec{n} = \begin{pmatrix} 0 \\ -10^4 \Big(1 - \cos\big(\frac{x_1 + 25}{50} \pi\big)\Big) \end{pmatrix} \qquad\text{on } (0,t_f) \times \Set{ \vec{x} \in \partial\mathcal{D} \given x_2 = 100 },
\end{equation*}
whereas no-flow condition and fixed displacement are given at the bottom
\begin{equation*}
    \kappa \nabla p \cdot \vec{n} = 0 \qquad\text{and}\qquad
    \vec{\eta} = \vec{0} \qquad\text{on } (0,t_f) \times \Set{ \vec{x} \in \partial\mathcal{D} \given x_2 = 0 }.
\end{equation*}
Furthermore, periodic boundary conditions are prescribed on the left and right boundaries at $x_1 \in \Set{ 0, 100 }$.
Setting the (dimensionless) Young's modulus $E = 10^5$ and Poisson's ratio $\nu = 0.9$ (with plane-stress assumption), we choose the physical parameters
\begin{equation*}
    \lambda = \frac{E \nu}{1 - \nu^2} = \frac{9 \cdot 10^6}{19} ,\qquad
    \mu = \frac{E}{2(1 + \nu)} = \frac{5 \cdot 10^5}{19} ,\qquad
    c_0 = 10^{-6} ,\qquad
    \kappa = 10^{-9},
\end{equation*}
and consider the strongly coupled situation where $\alpha = 1$.
Note that the incompatibility of the initial and boundary conditions again yields a solution which is discontinuous at time $t = 0$.

\begin{figure}[tb]%
    \centering
    \includegraphics[width=0.29\textwidth]{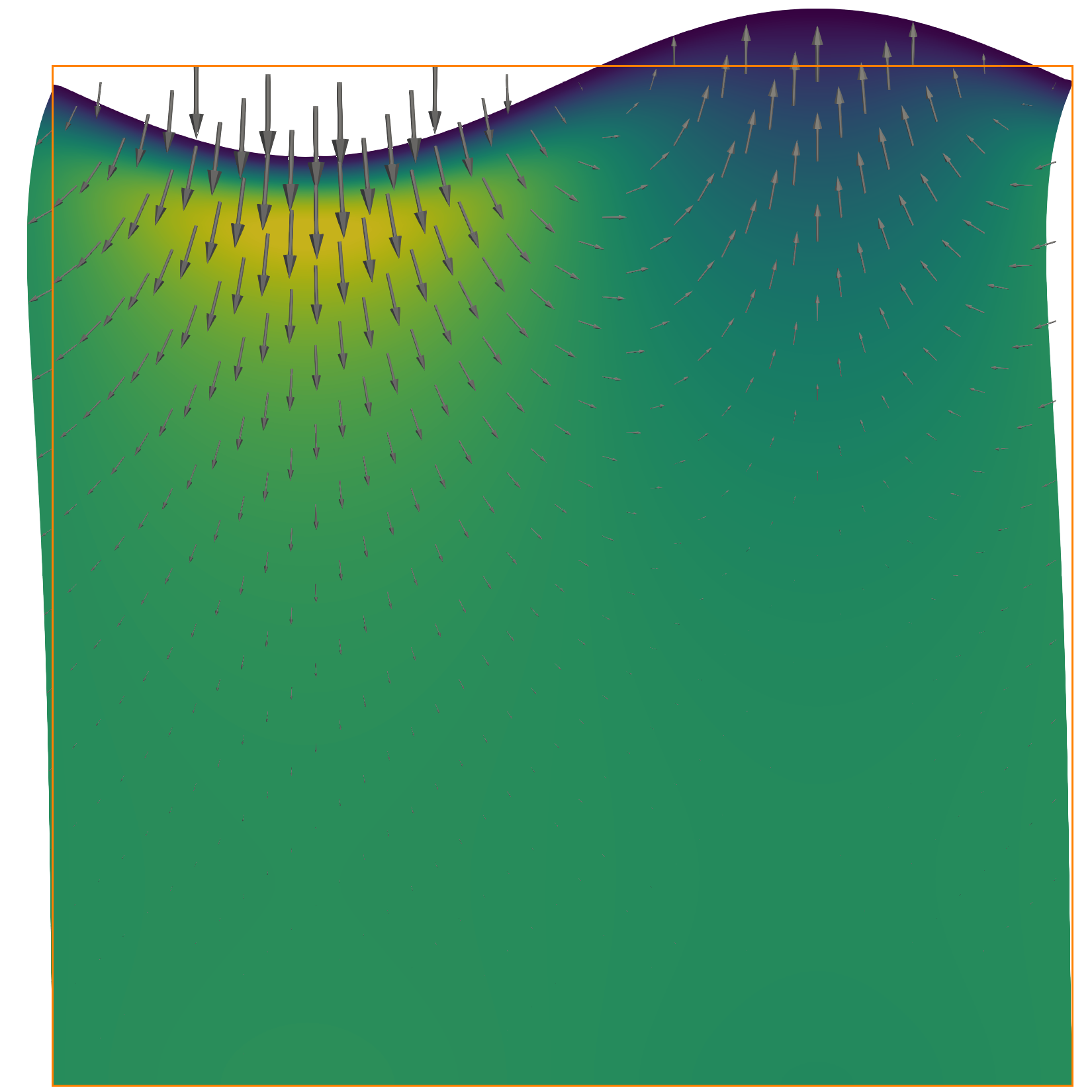}
    \includegraphics[width=0.29\textwidth]{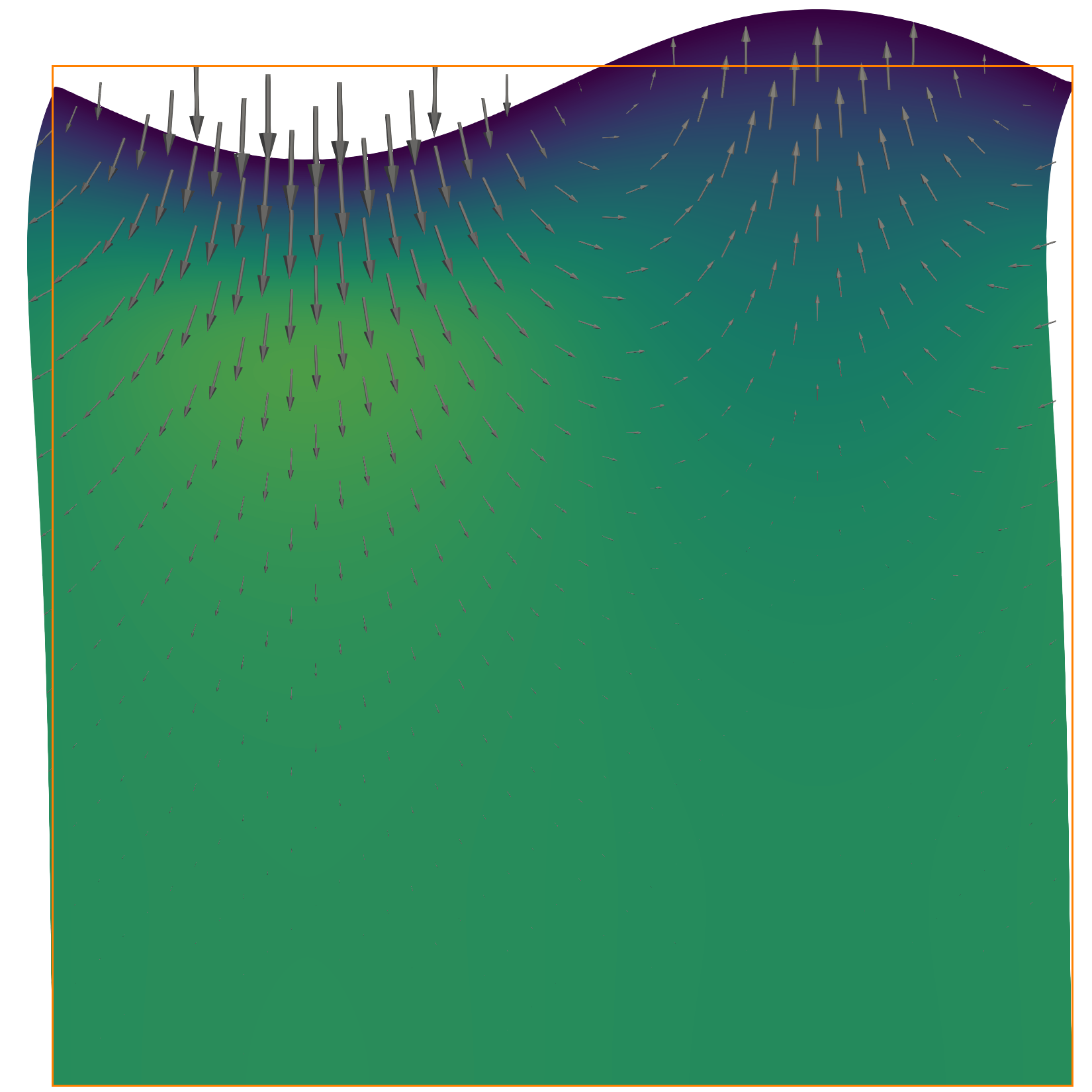}
    \includegraphics[width=0.29\textwidth]{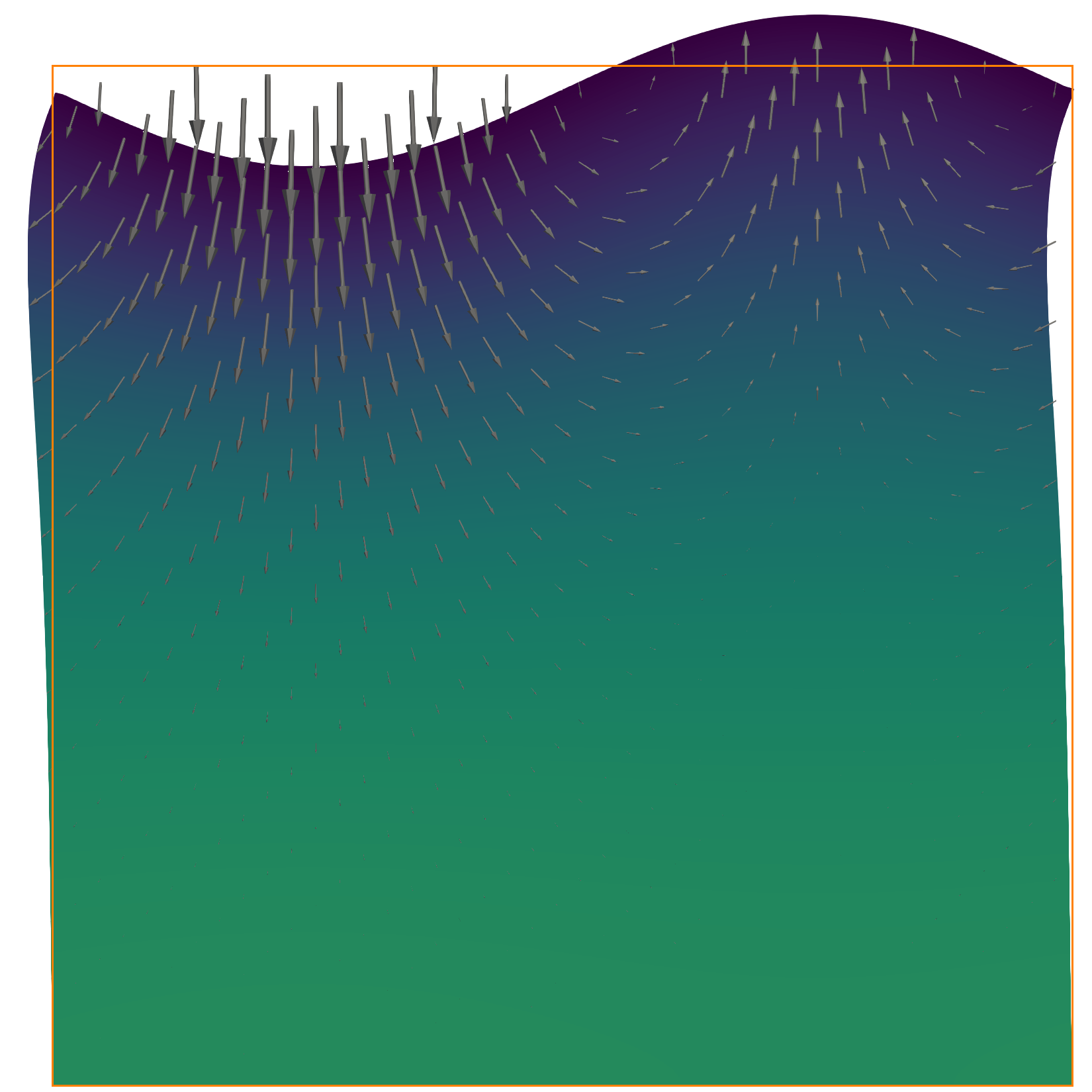}
    \includegraphics[width=0.07\textwidth]{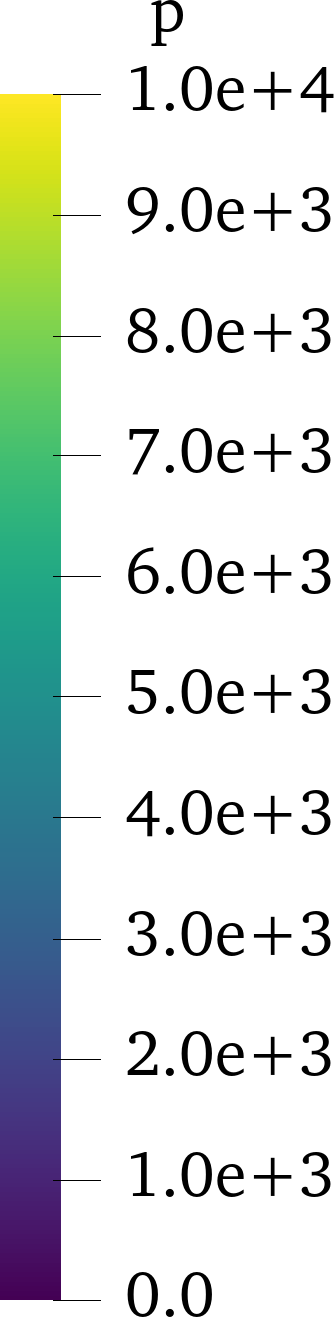}
    \caption{Numerical solutions at the times $t \in \Set[\big]{ 2 \cdot 10^4, 2 \cdot 10^5, 2 \cdot 10^6 }$ (left to right) for the 2D loading problem using the centered coupling scheme ($r = 1/2$) with MG(3,4) for $N_x = 128$.
        The pressure is represented in color, while the shown displacement is amplified by a factor of 3 and additionally represented by the arrows. The original outline is given by the solid orange line.}\label{fig:2D:solution}
\end{figure}

The numerical solution of the centered coupling scheme ($r = 1/2$) using MG(3,4) for $N_x = 128$ and $N_T = N_x^2 = 16384$ is depicted at the times $t \in \Set[\big]{ 2 \cdot 10^4, 2 \cdot 10^5, 2 \cdot 10^6 }$ in \cref{fig:2D:solution}.
Since the load is considerably larger on the left half, the subsidence is significant there, while the surface on the right is even lifted.
Furthermore, this leads to a horizontal expansion from the left half to the right one.
In the beginning, the fluid pressure is very large close to the left half of the surface, while it is low on the right half.
As expected, the pressure slowly vanishes over time resulting in a further subsidence of the whole surface.

\section{Conclusion}
\label{sec:conclusion}

In this work, we developed a novel coupling of LBMs to solve Biot's consolidation model in two dimensions.
To this end, we combined the classical LBM for reaction-diffusion equations \cite{Seta2013} with the recent pseudo-time LBM for quasi-static linear elasticity \cite{Boolakee2023}.
Therefore, the overall method uses only two distribution functions with (reduced) standard velocity sets for the solution fields, so that 17 unknowns per cell are required to solve for the 6 components of pressure, displacement and stress field at second-order consistency.
We addressed the coupling between the equations by a newly developed, centered update scheme which uses half of both explicit and semi-implicit contributions.

Numerical results obtained with the method of manufactured solutions confirm the analytical properties of the coupled methods.
This concerns the second-order convergence with respect to time-step size and lattice width.
Using a multi-grid method to accelerate the explicit pseudo-timestepping of the LBM for quasi-static elasticity, the required number of pseudo-time steps reduces from quadratic ($N_E \sim N_x^2$) to constant ($N_E \approx 3$ in the presented examples).
We have demonstrated that naively coupled LBMs lead to instabilities for strong coupling with Biot--Willis coefficient $\alpha$ close to 1, as both the fully explicit scheme and the semi-implicit scheme yield diverging solutions.
However, the centered coupling scheme is stable and accurate for all Biot--Willis coefficients.
Moreover, the presented numerical results for Terzaghi's (quasi-1D) consolidation problem and for a two-dimensional expansion thereof, both incorporate an instantaneous response due to external loading, and hence demonstrate that the coupled LBM is even applicable to discontinuous problems.

The method has several advantages, as it allows for a simple and scalable parallelization and is second-order consistent in space and time for all three solution fields.
The main drawback of the naive implementation lies in its pseudo-timestepping for quasi-static elasticity.
However, the presented multi-grid method eliminates this shortcoming since only few multi-grid steps are necessary.
Note that an adaptive number of multi-grid steps depending on the residual might reduce the computational cost even further.
Moreover, advanced multi-grid cycles, such as W- and F-cycle, as well as a full multi-grid approach are left for future research.
Considering complex geometries and the resulting implementation of boundary conditions in LBMs, the application of algebraic rather than geometric multi-grid might be preferable.
Lastly, also the Darcy flow subsystem can be accelerated analogously to allow for arbitrary time-step sizes.

In view of the fact, that this is the initial application of LBM to poroelasticity, there are ample opportunities for further development.
Besides possible extensions to three-dimensional and nonlinear poroelastic equations, different LBMs might be used to solve the flow and/or elasticity sub-problems.
In particular, solving the flow equations by a MRT scheme like in \cite{Ginzburg2005} could allow for anisotropic permeability, while accounting for elastodynamics solved by LBMs as e.g. in \cite{Escande2020,Boolakee2025} would enable the study of seismic waves in aquifers and reservoirs.
Furthermore, the limiting regimes of incompressible fluid and/or solid require appropriate treatment.
Finally, Biot's consolidation theory strongly simplifies real-world behavior, so that extensions towards nonlinear poro-(visco-)elastoplasticity including rate-dependency, as well as coupling to further processes, such as thermal exchange or chemical reactions, would open the methods applicability to many more contemporary applications.

\section*{CRediT authorship contribution statement}

\textbf{Stephan B. Lunowa:}
    Conceptualization,
    Data curation,
    Formal analysis,
    Investigation,
    Methodology,
    Software,
    Validation,
    Visualization,
    Writing – original draft,
    Writing – review and editing.
\textbf{Barbara Wohlmuth:}
    Conceptualization,
    Funding acquisition,
    Methodology,
    Resources,
    Supervision,
    Writing – review and editing.

\section*{Acknowledgments}

The authors thank Laura De Lorenzis, Oliver Boolakee and Martin Geier for the discussions about multi-grid methods for LBM.

\section*{Declaration of competing interest}

The authors declare that they have no known competing financial interests or personal relationships that could have appeared to influence the work reported in this paper.

\section*{Data availability}

The python implementation and scripts required to reproduce the numerical results of this article are available in the following GitHub repository: \url{https://github.com/s-lunowa/xlb-biot}.

\bibliographystyle{elsarticle-num}
\bibliography{literature}

\end{document}